\documentclass[10pt, twocolumn]{IEEEtran}

\usepackage{xspace,amsmath,amssymb,amsfonts,epsfig,syntonly}
\usepackage{cite,bm,color,url,textcomp}
\usepackage{epstopdf}
\usepackage{empheq}
\usepackage{flushend}

\newtheorem{remark}{\hspace{-11pt}\bf Remark}

\newcommand*\widefbox[1]{\fbox{\hspace{1em}#1\hspace{1em}}}

\DeclareMathOperator*{\argmin}{arg\,min}

\long\def\symbolfootnote[#1]#2{\begingroup
\def\thefootnote{\fnsymbol{footnote}}
\footnote[#1]{#2}\endgroup}

\psfull

\allowdisplaybreaks[4]

\begin{document}
\title{Robust Smart-Grid Powered \\ Cooperative Multipoint Systems}

\author{Xin Wang,~\IEEEmembership{Senior Member, IEEE}, Yu Zhang,~\IEEEmembership{Student Member, IEEE},\\
Georgios B. Giannakis,~\IEEEmembership{Fellow, IEEE}, and Shuyan Hu
\thanks {Work in this paper was supported by the the China Recruitment Program of Global Young Experts, the Program for New Century Excellent Talents in University, the Innovation Program of Shanghai Municipal Education Commission, the National Science and Technology Major Project of the Ministry of Science and Technology of China under Grant No.2012ZX03001013; US NSF grants CCF-1423316, CCF-1442686, ECCS-1202135,
and the Initiative for Renewable Energy \& the Environment (IREE) grant RL-0010-13.}
\thanks{X. Wang and S. Hu are with the Key Laboratory for Information Science of Electromagnetic Waves (MoE), Department of Communication Science and Engineering,
Fudan University, 220 Han Dan Road, Shanghai, China. Emails:~\{xwang11,\,syhu14\}fudan.edu.cn.

Y. Zhang and G. B. Giannakis are with the Department of Electrical and Computer Engineering and the Digital Technology Center,
University of Minnesota, Minneapolis, MN 55455, USA. Emails:~\{zhan1220,\,georgios\}@umn.edu.}
}

\markboth{} {}

\maketitle

\setcounter{page}{1}
\begin{abstract}
A framework is introduced to integrate renewable energy sources (RES) and dynamic pricing capabilities of the smart grid into beamforming designs for coordinated multi-point (CoMP) downlink communication systems. To this end, novel models are put forth to account for harvesting, storage of nondispatchable RES, time-varying energy pricing, as well as stochastic wireless channels. Building on these models, robust energy management and transmit-beamforming designs are developed to minimize the worst-case energy cost subject to the worst-case user QoS guarantees for the CoMP downlink. Leveraging pertinent tools, this task is formulated as a convex problem. A Lagrange dual based subgradient iteration is then employed to find the desired optimal energy-management strategy and transmit-beamforming vectors. Numerical results are provided to demonstrate the merits of the proposed robust designs.
\end{abstract}


\section*{Nomenclature}

\addcontentsline{toc}{section}{Nomenclature}

\subsection{Indices, numbers, and sets}

\begin{IEEEdescription}[\IEEEusemathlabelsep\IEEEsetlabelwidth{$M$, $m$}]

\item[$T$, $t$] Number and index of time slots.
\item[$I$, $i$] Number and index of distributed BSs.
\item[$K$, $k$] Number and index of mobile users.
\item[$S$, $s$] Number and index of sub-horizons.
\item[$j$] Iteration index of the dual subgradient ascent.
\item[$\mathcal{T}$] Set of the total scheduling horizon.
\item[$\mathcal{T}_{i,s}$] Sub-horizon $s$ for BS $i$.
\item[$\mathcal{K}$] Set of mobile users.
\item[$\mathcal{H}_k^t$] Channel uncertainty set of user $k$ in time slot $t$.
\item[$\mathcal{E}_i$] Uncertainty set of the energy harvested at BS $i$.
\item[$\mathcal{E}_i^{\textrm{p}}$, $\mathcal{E}_i^{\textrm{e}}$] Polyhedral and ellipsoidal cases of $\mathcal{E}_i$.
\item[$\ell$] Iteration index of the Bundle method.

\end{IEEEdescription}

\subsection{Constants}

\begin{IEEEdescription}[\IEEEusemathlabelsep\IEEEsetlabelwidth{$P_{D_n}^{\min,t}$, $P_{D_n}^{\max,t}$}]

\item[$\sigma_k^2$] Variance of channel additive noise.
\item[$\epsilon_k^t$] Radius of the channel uncertainty region for user $k$ in slot $t$.
\item[$\gamma_k$] Target signal-to-interference-plus-noise-ratio (SINR) of user $k$.
\item[$\underline{E}_{i}^t$, $\overline{E}_i^t$] Lower and upper bounds of the renewable energy at BS $i$ in slot $t$.
\item[$E^{\min}_{i,s}$, $E^{\max}_{i,s}$] Lower and upper bounds of the renewable energy at BS $i$ over the $s$th sub-horizon.
\item[$C_i^0$] Initial battery energy level of BS $i$.
\item[$C_i^{max}$] Battery capacity of BS $i$.
\item[$P_{b,i}^{\min}$, $P_{b,i}^{\max}$] Minimum and maximum (dis-)charging power at BS $i$.
\item[$\varpi_i$] Battery efficiency at BS $i$.
\item[$P_{c,i}$] Fixed energy consumption at BS $i$.
\item[$\xi$] Power amplifier efficiency.
\item[$P_{g,i}^{max}$] Upper limit of the total energy consumption for BS $i$.
\item[$\alpha^t$, $\beta^t$; $\psi^t$, $\phi^t$] Buying and selling energy prices; and functions thereof.

\end{IEEEdescription}

\subsection{Basic variables}
\begin{IEEEdescription}[\IEEEusemathlabelsep\IEEEsetlabelwidth{$\text{SINR}_k$}]

\item[$\mathbf{h}_{ik}^t$] Vector channel from BS $i$ to user $k$ in slot $t$.
\item[$\mathbf{h}_{k}^t$] Vector channel from all BSs to user $k$ in slot $t$.
\item[$\hat{\mathbf{h}}_k^t$] Estimated channel for user $k$ in slot $t$.
\item[$\mathbf{q}_k^t$] Vector signal transmitted to user $k$ in slot $t$.
\item[$s_k^t$] Information-bearing symbol of user $k$ in slot $t$.
\item[$y_k^t$] Received signal of user $k$ in slot $t$.
\item[$n_k^t$] Channel additive noise of user $k$ in slot $t$.
\item[$\mathbf{B}_i$] Selection matrix forming the $i$-th BS's transmit-beamforming vector.
\item[$\bm{\delta}_k^t$] Channel estimation error of user $k$ in slot $t$.
\item[$\mathbf{y}_{\ell}$, $\rho_{\ell}$] Proximal center and weight of the bundle method at iteration $\ell$.
\item[$\eta_{\ell}$] Predicted descent of the objective value for the bundle method at iteration $\ell$.
\item[$\mathbf{g}_{i,\ell}$] Subgradient of $\tilde{G}(\mathbf{p}_i)$ at iteration $\ell$.
\item[$r$] Energy buying-to-selling price ratio.
\item[$\hat{E}_i^t$] Expected renewable generation.
\item[$\tilde{E}_i^t$] Realizations of random renewable generation.
\item[$\kappa$, $U$] Scaling variable and uniform random variable relating to $\tilde{E}_i^t$.

\end{IEEEdescription}

\subsection{Decision variables}

\begin{IEEEdescription}[\IEEEusemathlabelsep\IEEEsetlabelwidth{$\text{SINR}_k$}]

\item[$\mathbf{w}_k^t$] Transmit beamforming vector from all BSs to user $k$ in slot $t$.
\item[$C_{i}^t$] Stored battery energy at BS $i$ at the beginning of slot $t$.
\item[$P_{b,i}^t$] Power delivered to or drawn from the batteries of BS $i$ in slot $t$.
\item[$P_{g,i}^t$] Total energy consumption for BS $i$ in slot $t$.
\item[$P_{x,i}^t$] Transmit power at BS $i$ in slot $t$.
\item[$E_i^t$] Energy harvested at BS $i$ in slot $t$.
\item[$\mathbf{e}_i$] Vector collecting $\{E_i^t\}_{t=1}^T$.
\item[$P_{i}^t$] Auxiliary variable relating $P_{g,i}^t$ and $P_{b,i}^t$.
\item[$\mathbf{p}_i$] Vector collecting $\{P_i^t\}_{t=1}^T$.
\item[$\mathbf{X}_k^t$] Matrix-lifting Beamformer for user $k$ in slot $t$.
\item[$\tau_k^t$] Auxiliary variable introduced by the S-procedure.
\item[$\mathbf{Z}$] Matrix collecting all primal variables.
\item[$\lambda_i^t$] Lagrange multiplier.
\item[$\mathbf{\Lambda}$] Matrix collecting all Lagrange multipliers.

\end{IEEEdescription}

\subsection{Functions}

\begin{IEEEdescription}[\IEEEusemathlabelsep\IEEEsetlabelwidth{$G(\cdot)$, $\tilde G(\cdot)$}]

\item[$\text{SINR}_k(\cdot)$] SINR of user $k$.
\item[$\widetilde{\text{SINR}}_k(\cdot)$] Worst-case SINR of user $k$.
\item[$G(\cdot,\cdot)$] Worst transaction cost across entire horizon.
\item[$G(\cdot)$, $\tilde G(\cdot)$] Modified worst-case transaction cost.
\item[$L(\mathbf{Z},\mathbf{\Lambda})$] Partial Lagrangian function.
\item[$D(\mathbf{\Lambda})$] Dual function.
\end{IEEEdescription}

\section{Introduction}

To accommodate the explosive demand for wireless services, cellular systems are evolving into what are termed heterogeneous networks (HetNets) consisting of distributed macro/micro/pico base stations (BSs) to cover overlapping areas of different sizes \cite{Hwa13}. Close proximity of many HetNet transmitters introduces severe inter-cell interference. For efficient interference management, coordinated multi-point processing (CoMP) has emerged as a promising technique for next-generation cellular networks such as LTE-Advanced \cite{Irm11}.

To fully exploit the potential of CoMP at affordable overhead, coordinated beamforming and/or clustered BS cooperation for downlink systems were investigated in \cite{Tol11, Dal10, Zha09, KimJainGG12, Ng10, Hong13}. Multiple BSs cooperate to beamform, and each user's data are only shared among a small number of BSs per cluster, thus greatly reducing the overall backhaul signalling cost.

The rapid development of small cells in HetNets has also driven the need for energy-efficient transmissions. Due to the growing number of BSs, the electricity bill has become a major part of the operational expenditure of cellular operators, and cellular networks contribute a considerable portion of the global ``carbon footprint'' \cite{Oh11}. These economic and ecological concerns advocate a ``green communication'' solution, where BSs in cellular networks are powered by the electricity grid \cite{Tol11, Dal10, Zha09, KimJainGG12, Oh11, Ng10, Hong13, Shi13, Rao13}. However, the current grid infrastructure is on the verge of a major paradigm shift, migrating from the traditional electricity grid to the so-termed ``smart grid'' equipped with a large number of advanced smart meters and
state-of-the-art communication and control links. To decrease greenhouse gas emissions, an important feature of future power systems is integration of renewable energy sources (RES). This leads to high penetration of distributed generators equipped with energy harvesting modules, which can crop energy from the environmental resources (e.g., solar and wind), and possibly trade the harvested energy with the main grid. In addition to distributed generation, distributed storage, and two-way energy trading associated with RES, demand-side management (DSM) including dynamic pricing and demand response, can further improve grid reliability and efficiency. Relying on pertinent tools, optimal energy management and scheduling with RES and/or DSM were proposed in \cite{Zha13, Liu10, Gua10, Gian13}.

To take advantage of the aforementioned smart grid capabilities in next-generation cellular systems, only a few recent works have considered the smart-grid powered CoMP transmissions \cite{Bu12, Xu13, Xu14}.
However, \cite{Bu12} only addressed dynamic pricing using a simplified smart grid level game, while \cite{Xu13} and~\cite{Xu14} assumed that the energy amounts harvested from RES are precisely available {\em a priori} (e.g., through forecasting), and the harvested energy cannot be stored at the BSs.
In addition, \cite{Xu14} assumed demand (or load) response based on different energy buying/selling prices across BSs without adapting power consumption to {\em time-varying} energy pricing.

The present paper deals with optimal energy management and transmit-beamforming designs for the smart-grid powered, cluster-based CoMP downlink with the clustering carried by the HetNet's central processor. Each BS has local RES and can perform two-way energy trading with the grid based on {\em dynamic} buying/selling prices. Different from \cite{Bu12, Xu13, Xu14}, we suppose that each BS has a local storage device, which can be charged to store the harvested (and even grid) energy and can be discharged to supply electricity if needed. To account for the stochastic and nondispatchable nature of both RES and wireless channels, we assume that the actual harvested energy amounts and the wireless channel states are {\em unknown} and {\em time-varying}, yet lie in some known uncertainty regions. Building on realistic models, we develop robust energy management and transmit-beamforming designs that minimize the worst-case energy cost subject to the worst-case user QoS guarantees for the CoMP downlink. Leveraging a novel formulation accounting for the worst-case transaction cost with two-way energy trading, as well as the S-procedure in robust beamforming designs, we show how to (re-)formulate the said task as a convex problem. Strong duality of the latter is then utilized to develop a Lagrange dual based subgradient solver. It is shown that the resultant algorithm is guaranteed to find the desired optimal energy-management strategy and transmit-beamforming vectors, and could also facilitate distributed implementations among the BSs.


The rest of the paper is organized as follows. The system models are described in Section II. The proposed approach to robust energy management and transmit-beamforming designs for the CoMP downlink is developed in Section III. Numerical results are provided in Section IV to demonstrate the proposed scheme. We conclude the paper in Section V.

\noindent {\em Notation}: Boldface fonts denote vectors or matrices; $\mathbb{C}^{K\times M}$ is the $K$-by-$M$ dimensional complex space; $(\cdot)'$ denotes transpose, and $(\cdot)^H$ conjugate transpose; $\text{tr}(\mathbf{A})$ and $\text{rank}(\mathbf{A})$ the trace and rank operators for matrix $\mathbf{A}$, respectively;
$\text{diag}(a_1, \ldots, a_M)$ denotes a diagonal matrix with $a_1, \ldots, a_M$ on the diagonal;
$| \cdot |$ represents the magnitude of a complex scalar; $\mathbf{A} \succeq \mathbf{0}$ means that a square matrix $\mathbf{A}$ is positive semi-definite.

\begin{figure*}
\centering \epsfig{file=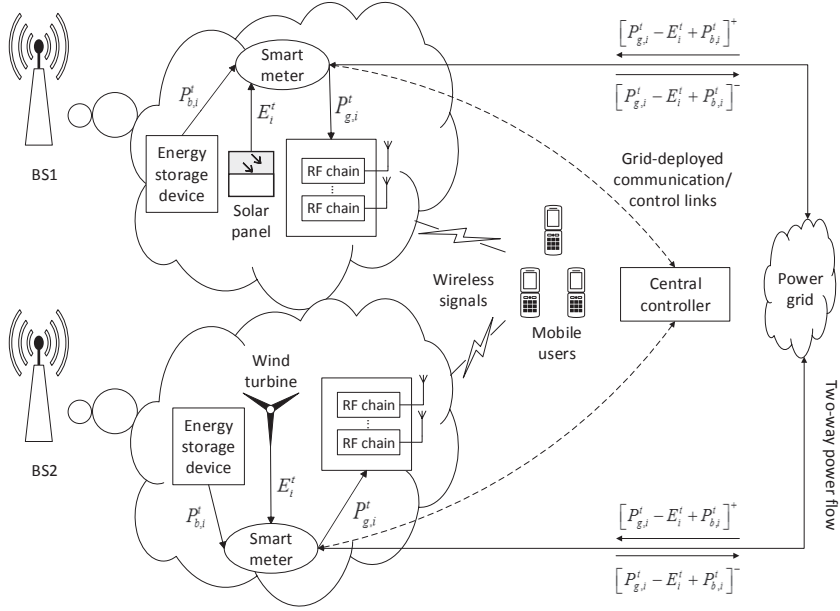, width=0.68\textwidth}
\caption{ A two-BS CoMP system powered by smart grids, where BSs with local renewable energy harvesting and storage devices implement two-way energy trading with the main grid. }
\label{F:model} 
\end{figure*}

\section{System Models}
\label{sec:model}

Consider a smart-grid powered cluster-based CoMP downlink. A set ${\cal I}:=\{1,\ldots, I\}$ of distributed BSs (e.g., macro/micro/pico BSs) provides service to a set ${\cal K}:=\{1,\ldots, K\}$ of mobile users; see Fig.~\ref{F:model}. Each BS has $M\geq 1$ transmit antennas and each user has a single receive antenna. Each BS is equipped with one or more energy harvesting devices (solar panels and/or wind turbines), and is also connected to the power grid with a two-way energy trading facility. Different from \cite{Bu12, Xu13, Xu14}, each BS has an energy storage device (possibly consisting of several large-capacity batteries) so that it does not have to consume or sell all the harvested energy on the spot, but can save it for later use.

For each CoMP cluster, there is a low-latency backhaul network connecting the set of BSs to a central controller \cite{Zha09}, which coordinates energy trading as well as cooperative communication. The central entity can collect both the communication data (transmit messages) from each of the BSs through the cellular backhaul links, as well as the energy information (energy buying/selling prices) from these BSs via smart meters installed at BSs, and the grid-deployed communication/control links connecting them.

Assume slot-based transmissions from the BSs to the users. While the actual harvested energy amounts and wireless channels cannot be accurately predicted, uncertainty regions for the wireless channels and renewable energy arrivals can be obtained, based on historical measurements and/or forecasting techniques. The slot duration is selected equal to the minimum time interval between the changes of the (channel or energy) uncertainty regions. Consider a finite scheduling horizon consisting of $T$ slots, indexed by the set ${\cal T}:= \{1, \ldots , T\}$. For convenience, the slot duration is normalized to unity unless otherwise specified; thus, the terms ``energy'' and ``power'' will hereafter be used interchangeably throughout the paper.

\subsection{Downlink CoMP Transmission Model}

Per slot $t$, let $\mathbf{h}_{ik}^t \in \mathbb{C}^{M \times 1}$ denote the vector channel from BS $i$ to user $k$, $\forall i \in {\cal I}$, $\forall k \in {\cal K}$; and let $\mathbf{h}_k^t := [{\mathbf{h}_{1k}^t}', \ldots, {\mathbf{h}_{Ik}^t}']'$ collect the channel vectors from all BSs to user $k$. With linear transmit beamforming performed across BSs, the vector signal transmitted  to user $k$ is
\[
    \mathbf{q}_k^t =\mathbf{w}_k^t s_k^t, \quad \forall k
\]
where $s_k^t$ denotes the information-bearing symbol, and $\mathbf{w}_k^t \in \mathbb{C}^{MI \times 1}$ denotes the beamforming vector across the BSs for user $k$. For convenience, $s_k^t$ is assumed to be a complex random variable with zero mean and unit variance. The received vector at user $k$ is therefore
\begin{equation}\label{eq.yk}
y_k^t ={\mathbf{h}_{k}^t}^H \mathbf{q}_k^t + \sum_{l\neq k} {\mathbf{h}_{k}^t}^H \mathbf{q}_l^t + n_k^t
\end{equation}
where ${\mathbf{h}_{k}^t}^H \mathbf{q}_k^t$ is the desired signal for user $k$, $\sum_{l\neq k} {\mathbf{h}_{k}^t}^H \mathbf{q}_l^t$ is the inter-user interference from the same cluster, and $n_k^t$ denotes additive noise, which may also include the downlink interference from other BSs outside this cluster. It is assumed that $n_k^t$ is a circularly symmetric complex Gaussian (CSCG) random variable with zero mean and variance $\sigma_k^2$.

The signal-to-interference-plus-noise-ratio (SINR) at user $k$ can be expressed as
\begin{equation}\label{sinr}
    \text{SINR}_k(\{\mathbf{w}_k^t\})= \frac{|{\mathbf{h}_{k}^t}^H \mathbf{w}_k^t|^2}{\sum_{l\neq k} (|{\mathbf{h}_{k}^t}^H \mathbf{w}_l^t|^2) + \sigma_k^2}~.
\end{equation}
The transmit power at each BS $i$ clearly is given by
\begin{equation}
    P_{x,i}^t = \sum_{k \in {\cal K}} {\mathbf{w}_k^t}^H \mathbf{B}_i \mathbf{w}_k^t
\end{equation}
where the matrix $\mathbf{B}_i\in \mathbb{R}^{MI\times MI}$
\begin{equation*}
\mathbf{B}_i:=\text{diag}\left(\underbrace{0, \ldots, 0}_{(i-1)M}, \underbrace{1, \ldots, 1}_{M}, \underbrace{0, \ldots, 0}_{(I-i)M}\right)
\end{equation*}
selects the corresponding rows out of $\{\mathbf{w}_k^t\}_{k\in \mathcal{K}}$ to form the $i$-th BS's transmit-beamforming vector.

The channel state information $\mathbf{h}_k^t$ is seldom precisely available {\it a priori} in practice. Relying on past channel measurements and/or reliable channel predictions, we adopt the following additive error model: $\mathbf{h}_k^t = \hat{\mathbf{h}}_k^t + {\mathbf{\delta}}_k^t$, where $\hat{\mathbf{h}}_k^t$ is the predicted channel. The uncertainty of this estimate is bounded by a spherical region \cite{Vucic2009}:
\begin{equation}\label{eq.chUncerty}
  \mathcal{H}_k^t:=\left\{\hat{\mathbf{h}}_k^t + {\mathbf{\delta}}_k^t \ | \  \left\|{\mathbf{\delta}}_k^t \right\| \leq \epsilon_k^t \right\}, \quad \forall k, t
\end{equation}
where $\epsilon_k^t>0$ specifies the radius of $\mathcal{H}_k^t$.
This leads to the worst-case SINR per user $k$ as [cf.~\eqref{sinr}]
\begin{equation}\label{eq.wcsnir}
  \widetilde{\text{SINR}}_k (\{\mathbf{w}_k^t\}) := \min_{\mathbf{h}_k^t \in\mathcal{H}_k^t} \frac{|{\mathbf{h}_{k}^t}^H \mathbf{w}_k^t|^2}{\sum_{l\neq k} (|{\mathbf{h}_{k}^t}^H \mathbf{w}_l^t|^2) + \sigma_k^2}~.
\end{equation}
To guarantee QoS per slot, it is required that
\begin{equation}\label{eq.sinr}
    \widetilde{\text{SINR}}_k (\{\mathbf{w}_k^t\}) \geq \gamma_k, \quad \forall k
\end{equation}
where $\gamma_k$ denotes the target SINR value per user $k$.

\subsection{Energy Harvesting Model}

Let $E_i^t$ denote the energy harvested during the last slot that is available at the beginning of slot $t$ at each BS $i \in {\cal I}$; and let $\mathbf{e}_i := [E_i^1 ,\ldots, E_i^T]'$. Due to the unpredictable and intermittent nature of RES, $\mathbf{e}_i$ is unknown a priori.
In general, uncertain quantities can be modeled by postulating either an underlying probability distribution
or an uncertainty region. Probability distributions (possibly mixed discrete/continuous) of the RES generation are seldom available in practice.
Although (non-)parametric approaches can be used to learn the distributions, the processes can be very complicated due to
the spatio-temporal correlations incurred by various meteorological factors.
On the other hand, the approach of postulating an uncertainty region provides the decision maker range forecasts instead of point forecasts,
which are essentially distribution-free.
Suppose that $\mathbf{e}_i$ lies in an uncertainty set ${\cal E}_i$, which can be obtained from historical measurements and/or fine forecast techniques.
From the perspective of computational tractability, two practical paradigms for ${\cal E}_i$ are considered here.
\begin{itemize}
\item[i)] The first model amounts to a polyhedral set \cite{Zha13}:
\begin{align}
{\cal E}_i^{\textrm{p}} := \Bigg\{\mathbf{e}_i\;|\; &\underline{E}_i^t \leq E_i^t \leq \overline{E}_i^t, \notag \\
&\hspace{-1cm} E_{i,s}^{\min} \leq \sum_{t \in {\cal T}_{i,s}} E_i^t \leq E_{i,s}^{\max}, \; {\cal T}=\bigcup_{s=1}^S {\cal T}_{i,s}\Bigg\} \label{eq.Ei}
\end{align}
where $\underline{E}_i^t$ ($\overline{E}_i^t$) denotes a lower (upper) bound on $E_i^t$; time horizon ${\cal T}$ is partitioned into consecutive but non-overlapping sub-horizons ${\cal T}_{i,s}$, $s = 1, \ldots, S$. Each sub-horizon can consist of multiple time slots, and the total energy harvested by BS $i$ over the $s$th sub-horizon is bounded by $E_{i,s}^{\min}$ and $E_{i,s}^{\max}$.

\item[ii)] The second model relies on an ellipsoidal uncertainty set (see e.g.~\cite{ChenDG13})
\begin{align}
{\cal E}_i^{\textrm{e}} := \left\{\mathbf{e}_i = \hat{\mathbf{e}}_i + \bm{\varsigma}_i \;
|\; \bm{\varsigma}_i'\bm{\Sigma}^{-1}\bm{\varsigma}_i \leq 1\right\} \label{eq.Ei2}
\end{align}
where $\hat{\mathbf{e}}_i := [\hat{E}_i^1 ,\ldots, \hat{E}_i^T]'$ denotes the nominal energy harvested at BS $i$,
which can be the forecasted energy, or simply its expected value. Vector $\bm{\varsigma}_i$ is the corresponding error in forecasting.
The known matrix $\bm{\Sigma}\succ \mathbf{0}$ quantifies the shape of the ellipsoid ${\cal E}_i^{\textrm{e}}$,
and hence determines the quality of the forecast.
\end{itemize}

\begin{remark}\textit{(Spatio-temporally correlated energy harvesting models)}.
The aforementioned two practical models capture RES uncertainty across the
scheduling (sub-)horizons per BS. The parameters required for constructing the sets
$\{{\cal E}_i^{\textrm{p}},{\cal E}_i^{\textrm{e}}\}$ can be obtained offline via statistical learning techniques using historical data.
In general, green energy harvested at different BSs can be spatially correlated if some BSs are geographically close.
In this case, joint spatio-temporal uncertainty models can be postulated whenever the underlying correlations are known {\it a priori}; see details in~\cite{Zha13}.
In general, a refined uncertainty model quantifying the actual harvested energy in a smaller region with a higher confidence level can be less conservative 
in the robust optimization formulation. However, the complexity of solving the resulting optimization problems directly depends on the choice of the uncertainty set.
\end{remark}

Note that the coherence time of RES arrivals can be much longer than that of wireless channels in practice \cite{Xu13, Xu15}. Yet, coherence times corresponding to the uncertainty regions of wireless channels can be much larger than those of the channel itself. In Section II-A, we implicitly assume that channel $\mathbf{h}_k^t$ remains unchanged per slot. However, $\widetilde{\text{SINR}}_k$ in (\ref{eq.wcsnir}) can be easily redefined as the worst-case SINR over multiple channel coherence times with the same uncertainty region (and possibly different channel realizations). This way the issue of different time scales becomes less critical. In addition, our models in fact accommodate cases where uncertainty regions for RES arrivals remain unchanged over multiple time slots. The proposed approach presented next readily applies to obtain the robust ahead-of-time schedule in this setup.

\subsection{Energy Storage Model}

Let $C_i^0$ denote the initial energy, and $C_i^t$ the amount of stored energy in the batteries of the $i$-th BS at the beginning of time slot $t$. With $C_i^{\max}$ bounding the capacity of batteries, it is clear that $0 \leq C_i^t \leq C_i^{\max}$, $\forall i \in {\cal I}$. Let $P_{b,i}^t$ denote the power delivered to or drawn from the batteries at slot $t$, which amounts to either charging ($P_{b,i}^t>0$) or discharging ($P_{b,i}^t<0$). Hence, the stored energy obeys the dynamic equation
\begin{equation}
    C_i^t = C_i^{t-1}+P_{b,i}^t, \quad t \in {\cal T}, \; i \in {\cal I}.
\end{equation}
The amount of power (dis-)charged is also bounded by
\begin{equation}
    \begin{split}
    P_{b,i}^{\min} \leq P_{b,i}^t \leq P_{b,i}^{\max}\\
    -\varpi_i C_i^{t-1} \leq P_{b,i}^t
    \end{split}
\end{equation}
where $P_{b,i}^{\min} <0$ and $P_{b,i}^{\max}>0$, while $\varpi_i \in (0,1]$ is the battery efficiency at BS $i$. The constraint $-\varpi_i C_i^{t-1} \leq P_{b,i}^t$ means that at most a fraction $\varpi_i$ of the stored energy $C_i^{t-1}$ is available for discharge.

\subsection{Energy Cost Model}

For the $i$-th BS per slot $t$, the total energy consumption $P_{g,i}^t$ includes the transmission-related power $P_{x,i}^t$, and the rest that is due to other components such as air conditioning, data processor, and circuits, which can be generally modeled as a constant power, $P_{c,i} > 0$ \cite{Xu14, Shi13}; namely,
\[
    P_{g,i}^t=P_{c,i}+P_{x,i}^t/\xi=P_{c,i}+\sum_{k \in {\cal K}} {\mathbf{w}_k^t}^H \mathbf{B}_i \mathbf{w}_k^t/\xi
\]
where $\xi>0$ denotes the power amplifier efficiency. For notational convenience, we absorb $\xi$ into $\mathbf{B}_i$ by redefining $\mathbf{B}_i := \mathbf{B}_i/\xi$; and further assume that $P_{g,i}^t$ is bounded by $P_{g,i}^{\max}$.

In addition to the harvested RES, the power grid can also supply the needed $P_{g,i}^t$ per BS $i$. With a two-way energy trading facility, the BS can also sell its surplus energy to the grid at a fair price in order to reduce operational costs. Given the required energy $P_{g,i}^t$, the harvested energy $E_i^t$, and the battery charging energy $P_{b,i}^t$, the shortage energy that needs to be purchased from the grid for BS $i$ is clearly $[P_{g,i}^t-E_i^t+P_{b,i}^t]^+$; or, the surplus energy (when the harvested energy is abundant) that can be sold to the grid is $[P_{g,i}^t-E_i^t+P_{b,i}^t]^-$, where   $[a]^+:= \max\{a, 0\}$, and $[a]^-  := \max\{-a, 0\}$. Note that both the shortage energy and surplus energy are non-negative, and we have at most one of them be positive at any time $t$ for \mbox{BS $i$}.

Suppose that the energy can be purchased from the grid at price $\alpha^t$, while the energy is sold to the grid at price $\beta^t$ per slot $t$. Notwithstanding, we shall always set $\alpha^t > \beta^t$, $\forall $, to avoid meaningless buy-and-sell activities of the BSs for profit. Assuming that the prices $\alpha^t$ and $\beta^t$ are known {\em a priori}, the \emph{worst-case transaction cost} for BS $i$ for the whole scheduling horizon is therefore given by
\begin{align}
G(\{P_{g,i}^t\},\{P_{b,i}^t\}):= \max_{\mathbf{e}_i \in {\cal E}_i} \sum_{t=1}^T \Big(\alpha^t[P_{g,i}^t-E_i^t+P_{b,i}^t]^+ \notag \\
- \beta^t [P_{g,i}^t-E_i^t+P_{b,i}^t]^-\Big).
\end{align}

\section{Energy Management for CoMP Beamforming}

Based on the models of Section~\ref{sec:model}, we consider here robust energy management for transmit beamforming in a CoMP cluster. Knowing only the uncertainty regions of the wireless channels and renewable energy arrivals, the central controller per cluster performs an (e.g. hour-) ahead-of-time schedule to optimize cooperative transmit beamforming vectors $\{\mathbf{w}_k^t\}$ and battery charging energy $\{P_{b,i}^t\}$, in order to minimize the worst-case total cost $\sum_{i \in {\cal I}} G(\{P_{g,i}^t\},\{P_{b,i}^t\})$, while satisfying the QoS guarantees $\widetilde{\text{SINR}}_k (\{\mathbf{w}_k^t\}) \geq \gamma_k$, $\forall k$, over the scheduling horizon ${\cal T}$. For convenience, we introduce the auxiliary variables $P_i^t:=P_{g,i}^t+P_{b,i}^t$, and formulate the problem as
\begin{subequations}\label{eq.prob}
\begin{align}
& \mathop{\text{minimize}}\limits_{\{\mathbf{w}_k^t,P_{b,i}^t,C_i^t,P_i^t\}} \;\sum_{i \in {\cal I}} G(\{P_i^t\}) \label{eq.proba} \\
& \text{suject to:} \notag \\
& P_i^t= P_{c,i} + \sum_{k \in {\cal K}} {\mathbf{w}_k^t}^H \mathbf{B}_i \mathbf{w}_k^t+P_{b,i}^t, \quad \forall i \in {\cal I}, \; t\in {\cal T} \label{eq.probb}\\
& 0 \leq P_{c,i} + \sum_{k \in {\cal K}} {\mathbf{w}_k^t}^H \mathbf{B}_i \mathbf{w}_k^t \leq P_{g,i}^{\max}, \quad \forall i \in {\cal I}, \; t\in {\cal T} \label{eq.probc}\\
  & C_i^t = C_i^{t-1} +P_{b,i}^t, \quad \forall i \in {\cal I}, \; t\in {\cal T}\label{eq.probd}\\
  &0 \leq C_i^t \leq C_i^{\max}, \quad \forall i \in {\cal I}, \; t\in {\cal T}\label{eq.probe}\\
  & P_{b,i}^{\min} \leq P_{b,i}^t \leq P_{b,i}^{\max}, \quad \forall i \in {\cal I}, \; t\in {\cal T}\label{eq.probf}\\
  & -\varpi_i C_i^{t-1} \leq P_{b,i}^t, \quad \forall i \in {\cal I}, \; t\in {\cal T}\label{eq.probg}\\
  & \widetilde{\text{SINR}}_k (\{\mathbf{w}_k^t\}) \geq \gamma_k, \quad \forall k \in {\cal K}, \; t\in {\cal T}. \label{eq.probh}
\end{align}
\end{subequations}
Consider for simplicity that (\ref{eq.prob}) is feasible. The problem can in fact become infeasible if the SINR thresholds $\gamma_k$ are high and the wireless channel qualities are not good enough. In this case, the subsequent problems \eqref{eq.sub1} are always infeasible. Recall that the proposed scheme is proposed to determine the ahead-of-time beamformer and energy schedule (offline). Such an infeasibility, once detected, can naturally lead to an admission control policy, i.e., to a criterion for dropping users or SINR requirements that render the problem infeasible \cite{Wan07}.

Solving \eqref{eq.prob} can provide a robust solution for the smart-grid powered CoMP downlink with worst-case performance guarantees.
It is worth mentioning here that thanks to the worst-case cost $G(\{P_i^t\})$, randomness introduced due to the wireless fading
propagation and also due to the RES uncertainty has been eliminated; thus~\eqref{eq.prob} contains only deterministic variables.
Because of~\eqref{eq.probb},~\eqref{eq.probc}, and~\eqref{eq.probh}, the problem is nonconvex, which motivates the reformulation
pursued next.

\subsection{Convex Reformulation}

First, let us consider convexity issues of the objective function $G(\{P_i^t\})$.
Define $\psi^t := (\alpha^t -\beta^t)/2$ and $\phi^t := (\alpha^t +\beta^t)/2$, and then rewrite:
\[
    G(\{P_i^t\})= \max_{\mathbf{e}_i \in {\cal E}_i} \sum_{t=1}^T \left(\psi^t|P_i^t-E_i^t| + \phi^t  (P_i^t-E_i^t)\right)~.
\]
Since $\alpha^t >\beta^t \geq 0$, we have $\phi^t > \psi^t >0$. It is then clear that $\psi^t|P_i^t-E_i^t| + \phi^t  (P_i^t-E_i^t)$ is a convex function of $P_i^t$ for any given $E_i^t$. As a pointwise maximization of these convex functions, $G(\{P_i^t\})$ is also a convex function of $\{P_i^t\}$ (even when the set ${\cal E}_i$ is non-convex) \cite{Bertsekas09}.

Except for \eqref{eq.probb},~\eqref{eq.probc}, and~\eqref{eq.probh}, all other constraints are convex. We next rely on the popular semidefinite program (SDP) relaxation technique to convexify (\ref{eq.probh}). By the definitions of $\mathcal{H}_k^t$ and $\widetilde{\text{SINR}}_k(\{\mathbf{w}_k^t\})$, the constraint $\widetilde{\text{SINR}}_k (\{\mathbf{w}_k^t\}) \geq \gamma_k$ can be rewritten as:
\begin{equation}\label{robustSINRconstraint}
F_k(\bm{\delta}_k^t)\geq 0~ \text{for all}~ \bm{\delta}_k^t ~\text{such that}~ {\bm{\delta}_k^t}^H\bm{\delta}_k^t \leq (\epsilon_k^t)^2
\end{equation}
where
\begin{equation*}
   F_k(\bm{\delta}_k^t):= (\hat{\mathbf{h}}_k^t+\bm{\delta}_k^t)^H(\frac{\mathbf{w}_k^t{\mathbf{w}_k^t}^H}{\gamma_k^t}-\sum_{l\neq k}\mathbf{w}_l^t{\mathbf{w}_l^t}^H)(\hat{\mathbf{h}}_k^t+\bm{\delta}_k^t) -\sigma_k^2.
\end{equation*}
Using~\eqref{robustSINRconstraint} and upon applying the well-known S-procedure in robust optimization \cite{Polik2007}, the original problem (\ref{eq.prob}) can be reformulated as an SDP with rank constraints.

\emph{(S-procedure)} Let $\mathbf{A, B}$ be $n \times n$ Hermitian matrices, $\mathbf{c}\in \mathbb{C}^n$ and $d\in \mathbb{R}$ for which the interior-point condition holds;
that is, there exists an $\bar{\mathbf{x}}$ such that
$\bar{\mathbf{x}}^H\mathbf{B}\bar{\mathbf{x}}<1$; and the following are equivalent:
\begin{itemize}
\item [(\emph{i})]
$ \mathbf{x}^H\mathbf{A}\mathbf{x}+\mathbf{c}^H\mathbf{x}+\mathbf{x}^H\mathbf{c}+d\geq 0$ for all $\mathbf{x}\in \mathbb{C}^n$ such that $\mathbf{x}^H\mathbf{B}\mathbf{x}\leq 1$;
\item[(\emph{ii})]
There exists a $\nu \ge 0$ such that
\begin{equation*}
\quad \left( \begin{matrix}
\mathbf{A}+\nu\mathbf{B} & \mathbf{c} \\
\mathbf{c}^H & d-\nu
\end{matrix}
\right) \succeq \mathbf{0}~.
\end{equation*}
\end{itemize}

For $\mathbf{X}_k^t:=\mathbf{w}_k^t {\mathbf{w}_k^t}^H\in \mathbb{C}^{MI\times MI}$, it clearly holds that $\mathbf{X}_k^t\succeq \mathbf{0}$ and $\text{rank}(\mathbf{X}_k^t)=1$. Using the S-procedure,~\eqref{robustSINRconstraint} can be transformed to
 \begin{equation}\label{S-procedure-Cons}
 \begin{split}
 \mathbf{\Gamma}_k^t:=\left(\begin{matrix}
       \mathbf{Y}_k^t+\tau_k^t\mathbf{I} & \mathbf{Y}_k^t\hat{\mathbf{h}}_k^t \\
       \hat{\mathbf{h}}_k^{tH}\mathbf{Y}_k^{tH}   & \hat{\mathbf{h}}_k^{tH}\mathbf{Y}_k^t\hat{\mathbf{h}}_k^t-\sigma_k^2-\tau_k^t(\epsilon_k^t)^2 \end{matrix} \right) \succeq \mathbf{0}
\end{split}
\end{equation}
where $\tau_k^t\ge 0$ and
\begin{equation}
\mathbf{Y}_k^t:=\frac{1}{\gamma_k}\mathbf{X}_k^t-\sum_{l\neq k} \mathbf{X}_l^t.
\end{equation}

Introducing auxiliary variables $\tau_k^t$ and dropping the rank constraints $\text{rank}(\mathbf{X}_k^t)=1$, $\forall k,t$, we can then relax (\ref{eq.prob}) to:
\begin{subequations}
\label{eq.sdr}
\begin{empheq}[box=\widefbox]{align}
&\min_{\{\mathbf{X}_k^t,\tau_k^t, P_i^t,
P_{b,i}^t,C_i^t\}} \, \sum_{i \in {\cal I}} G(\{P_i^t\})  \label{eq.Reforma}\\
&\text{s. t.}~~\text{(\ref{eq.probd})--(\ref{eq.probg})} \notag \\
& P_i^t = P_{c,i}+\sum_{k \in {\cal K}} \text{tr}(\mathbf{B}_i\mathbf{X}_k^t) + P_{b,i}^t, \, \forall i, t \label{eq.Reformb} \\
& 0 \leq \sum_{k \in {\cal K}} \text{tr}(\mathbf{B}_i\mathbf{X}_k^t) \leq P_{g,i}^{\max}-P_{c,i}, \, \forall i, t \label{eq.Reformc} \\
& \mathbf{\Gamma}_k^t\succeq \mathbf{0}, \;\; \tau_k^t \geq 0, \;\; \mathbf{X}_k^t \succeq \mathbf{0}, \, \forall k,t. \label{eq.Reformd}
\end{empheq}
\end{subequations}
In addition to the linear constraints (\ref{eq.probd})--(\ref{eq.probg}), the quadratic power constraints \eqref{eq.probb} and \eqref{eq.probc} have now become linear, and SINR constraints \eqref{eq.probh} become a set of convex SDP constraints in (\ref{eq.sdr}). Since the objective function is convex, problem~\eqref{eq.sdr} is a convex program that can be tackled by a centralized solver; e.g., using the projected subgradient descent approach.
However, the feasible set~\eqref{eq.Reformb}--\eqref{eq.Reformd} is the intersection of a semi-definite cone and a polytope, for which the iterative projection is complicated. To reduce computational complexity and enhance resilience to failures, we next develop an efficient algorithm to solve~\eqref{eq.sdr} in a distributed fashion coordinated by different agents.

\subsection{Lagrange Dual Approach}

Since (\ref{eq.sdr}) is convex, a Lagrange dual approach can be developed to efficiently find its solution. Let $\lambda_i^t$, $\forall i,t$ denote the Lagrange multipliers associated with the constraints~\eqref{eq.Reformb}. With the convenient notation $\mathbf{Z}:=\{\mathbf{X}_k^t,\tau_k^t, P_{b,i}^t,C_i^t,P_i^t\}$ and $\mathbf{\Lambda}:=\{\lambda_i^t\}$, the partial Lagrangian function of (\ref{eq.sdr}) is
\begin{align}
&L(\mathbf{Z},\mathbf{\Lambda}) := \notag \\
&\sum_{i \in {\cal I}} \Bigl[G(\{P_i^t\})+ \sum_{t=1}^T\lambda_i^t\bigl(P_{c,i}+\sum_{k \in {\cal K}} \text{tr}(\mathbf{B}_i\mathbf{X}_k^t)+P_{b,i}^t-P_i^t\bigr)\Bigr].~\label{eq.Lam}
\end{align}
The Lagrange dual function is then given by
\begin{equation}
\begin{split}
D(\mathbf{\Lambda})& :=\min_{\mathbf{Z}}\; L(\mathbf{Z},\mathbf{\Lambda})\\
 \text{s. t.} & ~\eqref{eq.probd}-\eqref{eq.probg},\,\eqref{eq.Reformc}-\eqref{eq.Reformd}\label{eq.dual func}
\end{split}
\end{equation}
and the dual problem of~\eqref{eq.sdr} is:
\begin{equation}
\max_{\mathbf{\Lambda}}\; D(\mathbf{\Lambda}).\label{eq.dual prob}
\end{equation}

{\em Subgradient iterations}: Let $j$ denote the iteration index.
To obtain the optimal solution $\mathbf{\Lambda}^*$ to the dual problem \eqref{eq.dual prob},
we resort to the dual subgradient ascent method, which amounts to the following update
\begin{align}
\label{eq.subgradient}
\lambda_i^t(j+1) = \lambda_i^t(j)+ \mu(j) g_{\lambda_i^t}(j),\,\forall i,t
\end{align}
where $\mu(j)>0$ is an appropriate stepsize. The subgradient $\mathbf{g}(j):=[g_{\lambda_i^t}(j)\; \forall i,t]$ can then be expressed as
\begin{equation}\label{eq.subgMlpr}
g_{\lambda_i^t}(j)= P_{c,i}+\sum_{k \in {\cal K}} \text{tr}(\mathbf{B}_i\mathbf{X}_k^t(j))+P_{b,i}^t(j)-P_i^t(j),\, \forall i,t
\end{equation}
where $\mathbf{X}_k^t(j)$, $P_{b,i}^t(j)$, and $P_i^t(j)$ are given by
\begin{subequations}
\begin{align}
 \{\mathbf{X}_k^t(j)\}_{k=1}^K & \in \arg \;\min_{\{\mathbf{X}_k^t, \tau_k^t\}} \; \sum_{i \in {\cal I}}[\lambda_i^t(j) \sum_{k \in {\cal K}} \text{tr}(\mathbf{B}_i\mathbf{X}_k^t)] \nonumber\\
& \text{s. t.} ~~ \eqref{eq.Reformc}-\eqref{eq.Reformd},\;\; \forall k \label{eq.sub1}\\
\{P_{b,i}^t(j)\}_{t=1}^T & \in \arg \;\min_{\{P_{b,i}^t, C_{i,t}\}} \; \sum_{t=1}^T [\lambda_i^t(j)P_{b,i}^t] \nonumber\\
&\text{s. t.} ~~\eqref{eq.probd}-\eqref{eq.probg}, \;\; \forall t \label{eq.sub2}\\
\{P_i^t(j)\}_{t=1}^T & \in \arg \;\min_{\{P_i^t\}} \; [G(\{P_i^t\})-\sum_{t=1}^T\lambda_i^t(j)P_i^t]~. \label{eq.sub3}
\end{align}
\end{subequations}

{\em Solving the subproblems}: Subproblems~\eqref{eq.sub1} are standard SDPs per $t \in {\cal T}$; hence, $\{\mathbf{X}_k^t(j)\}_{k=1}^K$ for all $t$ can be efficiently solved by general interior-point methods \cite{Bertsekas09}.

The subproblems \eqref{eq.sub2} are simple linear programs (LPs) over $\{P_{b,i}^t, C_i^t\}_{t=1}^T$ per $i \in {\cal I}$; hence, $\{P_{b,i}^t(j)\}_{t=1}^T$, $\forall i$, can be obtained by existing efficient LP solvers.

Due to the convexity of $G(\{P_i^t\})$, the subproblems \eqref{eq.sub3} are convex per $i \in {\cal I}$. Yet, because of the non-differentiability of $G(\{P_i^t\})$ due to the absolute value operator and the maximization over $\mathbf{e}_i \in {\cal E}_i$, the problem is challenging to be handled by existing general solvers. For this reason, we resort to the proximal bundle method to obtain $\{P_i^t(j)\}_{t=1}^T$. Upon defining $\tilde{G}(\mathbf{p}_i):=G(\mathbf{p}_i)-\sum_{t=1}^T\lambda_i^t(j)P_i^t$,
where $\mathbf{p}_i := [P_i^1,\ldots,P_i^T]'$, the partial subgradient of $\tilde{G}(\mathbf{p}_i)$  with respect to $P_i^t$ can be obtained as
\begin{equation}\label{eq.pG}
    \frac{\partial \tilde{G}(\mathbf{p}_i)}{\partial P_i^t} =
    \begin{cases}
        \alpha_i^t - \lambda_i^t(j), & \text{if  } P_i^t\geq {E_i^t}^*\\
        \beta_i^t - \lambda_i^t(j), & \text{if  } P_i^t< {E_i^t}^*
    \end{cases}
\end{equation}
where $\mathbf{e}_i^*:=[{E_i^1}^*,\ldots, {E_i^T}^*]'$ for the given $\mathbf{p}_i$ is obtained as
\begin{equation}\label{subprob.maxE}
    \mathbf{e}_i^* \in \arg\; \max_{\mathbf{e}_i \in {\cal E}_i}\; \sum_{t=1}^T \left(\psi^t|P_i^t-E_i^t| + \phi^t  (P_i^t-E_i^t)\right).
\end{equation}
It can be readily checked that the objective function in~\eqref{subprob.maxE} is convex in~$\mathbf{e}_i$ under the condition
$\alpha^t > \beta^t,\, \forall t\in \mathcal{T}$. Hence, the globally optimal solution is
attainable at the extreme points of ${\cal E}_i$. Leveraging the special structure of ${\cal E}_i$,
we utilize an efficient vertex enumerating algorithms to obtain $\mathbf{e}_i^*$ directly,
as detailed next.

\subsection{Proximal Bundle Method}
Given the partial subgradient in~\eqref{eq.pG}, nonsmooth convex optimization algorithms
can be employed to solve the subproblem~\eqref{eq.sub3}. A state-of-the-art bundle method~\cite[Ch.~6]{Bertsekas09},~\cite{Kiwiel00}
will be developed here with guaranteed convergence to the optimal $\{P_i^t(j)\}_{t=1}^T$; see also~\cite{Zha13}.

Similar to cutting plane methods, the idea of bundle method is to approximate
the epigraph of an objective by the intersection of a number of
halfspaces. The latter are generated through the supporting hyperplanes, referred to as cuts,
by using the subgradients. Specifically, letting $\ell$ denote the iteration index of the bundle method,
the iterate $\mathbf{p}_{i,\ell+1}$ is obtained by minimizing a polyhedral
(piecewise linear) approximation of $\tilde{G}(\mathbf{p}_i)$ with a quadratic proximal regularizer
\begin{align}
\label{prob:bundle}
\mathbf{p}_{i,\ell+1} :=\argmin\limits_{\mathbf{p}_i \in \mathbb{R}^T} \left\{\hat{G}_{\ell}(\mathbf{p}_i)+\frac{\rho_{\ell}}{2}\|\mathbf{p}_i-\mathbf{y}_{\ell}\|^2\right\}
\end{align}
where
$\hat{G}_{\ell}(\mathbf{p}_i):=\max\{\tilde{G}(\mathbf{p}_{i,0})+\mathbf{g}_{0}^{\prime}(\mathbf{p}_i-\mathbf{p}_{i,0}),
\ldots,\tilde{G}(\mathbf{p}_{i,\ell})+\mathbf{g}_{\ell}^{\prime}(\mathbf{p}_i-\mathbf{p}_{i,\ell})\}$;
$\mathbf{g}_{i,\ell}$ is the subgradient of $\tilde{G}(\mathbf{p}_i)$
evaluated at the point $\mathbf{p}=\mathbf{p}_{i,\ell}$ [cf.~\eqref{eq.pG}];
and the proximity weight $\rho_{\ell}$ controls the stability of iterates.

Different from the proximal cutting plane method (CPM), the bundle method updates its
proximal center $\mathbf{y}_{\ell}$  according to a descent query
\begin{align}\label{eq.query}
\mathbf{y}_{\ell+1}  =
\left\{\begin{array}{cc}
\mathbf{p}_{i,\ell+1},  &\mbox{if}~\tilde{G}(\mathbf{y}_{\ell})-\tilde{G}(\mathbf{p}_{i,\ell+1}) \ge \theta\eta_{\ell} \\
\mathbf{y}_{\ell},  &\mbox{if}~\tilde{G}(\mathbf{y}_{\ell})-\tilde{G}(\mathbf{p}_{i,\ell+1}) < \theta\eta_{\ell}
\end{array}\right.
\end{align}
where $\theta \in (0,1)$ is a pre-selected constant,  and $\eta_{\ell} := \tilde{G}(\mathbf{y}_{\ell})-\left(\hat{G}_{\ell}(\mathbf{p}_{\ell+1})+\frac{\rho_{\ell}}{2}\|\mathbf{p}_{\ell+1}-\mathbf{y}_{\ell}\|^2\right)$ is the predicted descent of the objective in \eqref{prob:bundle}.  Essentially, if the actual descent amount $\tilde{G}(\mathbf{y}_{\ell})-\tilde{G}(\mathbf{p}_{\ell+1})$ is no less than a $\theta$ fraction of the predicted counterpart $\eta_{\ell} $,
then the iterate takes a ``serious'' step updating its proximal center $\mathbf{y}_{\ell+1}$  to the latest point $\mathbf{p}_{\ell+1}$; otherwise it is just a ``null'' step with the center unchanged. The intelligent query~\eqref{eq.query} enables the bundle method to find ``good'' proximal centers along the iterates, and hence it converges faster than the proximal CPM.
In addition, depending on whether a serious or a null step is taken, the proximity weight $\rho_{\ell}$ can be updated accordingly to further accelerate convergence~\cite{Kiwiel95};  that is,
\begin{align*}
\rho_{\ell+1}  =
\left\{\begin{array}{cc}
 \max(\rho_{\ell}/10, \rho_{\min}),\, & \mbox{if}~\tilde{G}(\mathbf{y}_{\ell})-\tilde{G}(\mathbf{p}_{i,\ell+1}) \ge \theta\eta_{\ell} \\
 \min(10\rho_{\ell}, \rho_{\max}),\,  &\mbox{if}~\tilde{G}(\mathbf{y}_{\ell})-\tilde{G}(\mathbf{p}_{i,\ell+1}) < \theta\eta_{\ell} .
\end{array}\right.
\end{align*}
The algorithm terminates if $\mathbf{y}_{\ell}=\mathbf{p}_{i,\ell+1}$, while finite termination
is achievable if both the objective and the feasible set are polyhedral~\cite[Sec.~6.5.3]{Bertsekas09}.

Now, to complete the proximal bundle method for solving \eqref{eq.sub3}, we only need solve problem \eqref{prob:bundle}. Using an auxiliary variable $u$, \eqref{prob:bundle} can be re-written as
\begin{subequations}
\label{prob:QPprimal}
\begin{align}
&\min\limits_{\mathbf{p}_i, u}~\quad u+\frac{\rho_{\ell}}{2}\|\mathbf{p}_i-\mathbf{y}_{\ell}\|^2 \\
&\text{s.t.}~\,\, \tilde{G}(\mathbf{p}_{i,n})+\mathbf{g}_{i,n}^{\prime}(\mathbf{p}_i-\mathbf{p}_{i,n}) \le u,~n = 0,1,\ldots,\ell.
\end{align}
\end{subequations}
Introducing multipliers $\bm{\xi} \in \mathbb{R}_{+}^{\ell+1}$ and
letting $\mathbf{1}$ denote the all-ones vector, the Lagrangian function is given as
\begin{align}
\label{QPLagrang}
\mathcal{L}(u, \mathbf{p}_i, \bm{\xi}) = (1-\mathbf{1}'\bm{\xi})u+\frac{\rho_{\ell}}{2}\|\mathbf{p}_i-\mathbf{y}_{\ell}\|^2 \nonumber \\
+\sum_{n=0}^{\ell}\xi_n\big(\tilde{G}(\mathbf{p}_{i,n})+\mathbf{g}_{i,n}^{\prime}(\mathbf{p}_i-\mathbf{p}_{i,n})\big).
\end{align}
The optimality condition $\nabla_{\mathbf{p}_i}\mathcal{L}(u, \mathbf{p}_i, \bm{\xi})=\mathbf{0}$
yields
\begin{align}\label{OptP}
\mathbf{p}_i^{*} =  \mathbf{y}_{\ell}-\frac{1}{\rho_{\ell}}\sum_{n=0}^{\ell}\xi_n\mathbf{g}_{i,n}.
\end{align}
Substituting \eqref{OptP} into \eqref{QPLagrang}, the dual of
\eqref{prob:QPprimal} is
\begin{subequations}
\label{QPdual}
\begin{align}
&\hspace{-4.0mm} \max\limits_{\bm{\xi}} \,  -\frac{1}{2\rho_{\ell}}\left\|\sum_{n=0}^{\ell}\xi_n\mathbf{g}_{i,n}\right\|^2 +
\sum_{n=0}^{\ell}\xi_n\big(\tilde{G}(\mathbf{p}_{i,n})+\mathbf{g}_{i,n}^{\prime}(\mathbf{y}_{\ell}-\mathbf{p}_{i,n})\big)\\
&\text{s.t.} ~\quad  \bm{\xi}\succeq \mathbf{0},~\mathbf{1}'\bm{\xi}= 1~.
\end{align}
\end{subequations}
It can be readily seen that \eqref{QPdual} is essentially a QP over the probability simplex in the dual space,  which can be solved very efficiently as in e.g.,~\cite{Bertsekas82}.

\subsection{Optimality and Distributed Implementation}

When a constant stepsize $\mu(j)=\mu$ is adopted, the subgradient iterations (\ref{eq.subgradient}) are guaranteed to converge to a neighborhood of the optimal $\mathbf{\Lambda}^*$ for the dual problem (\ref{eq.dual prob}) from any initial $\mathbf{\Lambda}(0)$. The size of the neighborhood is proportional to the stepsize $\mu$. If we adopt a sequence of non-summable diminishing stepsizes satisfying $\lim_{j \rightarrow \infty} \mu(j) =0$ and $\sum_{j=0}^{\infty} \mu(j) = \infty$, then the iterations (\ref{eq.subgradient}) will asymptotically converge to the exact $\mathbf{\Lambda}^*$ as $j \rightarrow \infty$ \cite{Bertsekas09}.

The objective function \eqref{eq.Reforma} is not strictly convex because it does not involve
all optimization variables. Hence, when it comes to primal convergence,
extra care is necessary~\cite{Xiao04,Gatsis12}.
Specifically, the optimal primal can be attained either by adding a strictly convex regularization term,
or, by utilizing the augmented Lagrangian.
Here, we will simply implement the~\emph{Ces\'{a}ro averaging} method~\cite{Kiwiel04}
to obtain the optimal power schedules.
With $\mu_{\text{sum}}^m:=\sum_{j=1}^m\mu(j)$, the running average is
\begin{align}\label{eq:cesaro}
\bar{\mathbf{Z}}^m := \frac{1}{\mu_{\text{sum}}^m}\sum_{j=1}^m\mu(j)\mathbf{Z}(j)
\end{align}
which can be efficiently computed in a recursive way
\begin{align}\label{eq:cesaro2}
\bar{\mathbf{Z}}^m := \frac{\mu(m)}{\mu_{\text{sum}}^m}\mathbf{Z}(m)+
 \frac{\mu_{\text{sum}}^{m-1}}{\mu_{\text{sum}}^m}\bar{\mathbf{Z}}^{m-1}~.
\end{align}
Note that if a constant stepsize $\mu(j)\equiv \mu$ is adopted,~\eqref{eq:cesaro} and
~\eqref{eq:cesaro2} boil down to the ordinary running average
 \begin{align*}
\bar{\mathbf{Z}}^m = \frac{1}{m}\sum_{j=1}^m\mathbf{Z}(j) =
\frac{1}{m}\mathbf{Z}(m) + \frac{m-1}{m}\bar{\mathbf{Z}}^{m-1}\,.
\end{align*}

If for the relaxed problem (\ref{eq.sdr}), the obtained solution satisfies the condition $\text{rank}({\mathbf{X}_k^t}^*)=1$, $\forall k,t$,
then it clearly yields the optimal beamforming vectors ${\mathbf{w}_k^t}^*$ as the (scaled) eigenvector with respect to the only positive eigenvalue of ${\mathbf{X}_k^t}^*$ for the original problem (\ref{eq.prob}). Fortunately, it was shown in \cite[Them. 1]{Song2012} that the S-procedure based SDP (\ref{eq.sdr}) always returns a rank-one optimal solution ${\mathbf{X}_k^t}^*$, $\forall k,t$, when the uncertainty bounds $\epsilon_k^t$ are sufficiently small. If $\epsilon_k$ is large, the existence of rank-one optimal solutions of (\ref{eq.sdr}) cannot be provably guaranteed. In this case, a randomized rounding strategy \cite{Tseng03} needs to be adopted to obtain vectors ${\mathbf{w}_k^t}^*$ from ${\mathbf{X}_k^t}^*$ to nicely approximate the solution of the original problem (\ref{eq.prob}). Even though no proof is available to ensure a rank-one solution when $\epsilon_k$ is large,
it has been extensively observed in simulations that the SDP relaxation always returns a rank-one optimal solution~\cite{Song2012}. This confirms the view that the optimal beamforming vectors for the original problem (\ref{eq.prob})
will be obtained by our approach with high probability.

The subgradient iterations can be run in a distributed fashion. Specifically, the central controller can maintain the Lagrange multipliers $\mathbf{\Lambda}(j)$ and broadcast them to all BSs via backhaul links. Given the current $\mathbf{\Lambda}(j)$, the central controller solves the subproblems (\ref{eq.sub1}) to obtain the beamforming vectors for all BSs. On the other hand, each BS solves its own subproblems (\ref{eq.sub2})--(\ref{eq.sub3}), which are decoupled across  BSs. The BSs send back to the central controller $\{P_{g,i}^t(j)\}_{t=1}^T$, $\{P_{b,i}^t(j)\}_{t=1}^T$, and $\{P_i^t(j)\}_{t=1}^T$, which are in turn used to update $\mathbf{\Lambda}(j+1)$ through the subgradient iterations (\ref{eq.subgradient}).

%
\section{Numerical Tests}
\label{sec:test}
%

\begin{table}[t]
\centering
\caption{Generating capacities, battery initial energy and capacity,   charging limits and efficiency.
}\label{tab:param}
    \begin{tabular}{  c | c  c  c  c  c  c  c }
    \hline
Unit &$P_{G_i}^{\min}$ &$P_{G_i}^{\max}$   &$C_i^{0}$ &$C_i^{\max}$  &$P_{b_i}^{\min}$     &$P_{b_i}^{\max}$   &$\varpi_i$ \\ \hline
1    & 0        & 50           & 5  	  & 30             &-10               & 10          & 0.95             \\
2    & 0        & 45          & 5       & 30           &-10                 & 10           & 0.95          \\
3    & 0        & 45           & 5  	  & 30             &-10               & 10          & 0.95             \\
4    & 0        & 45          & 5       & 30           &-10                 & 10           & 0.95          \\
5    & 0        & 50           & 5  	  & 30             &-10               & 10          & 0.95             \\
6    & 0        & 45          & 5       & 30           &-10                 & 10           & 0.95          \\
 \hline
    \end{tabular}
\end{table}

\begin{table}[t]
\renewcommand{\arraystretch}{1.2}
\centering
\caption{Limits of forecasted wind power and energy purchase prices}\label{tab:wind}
    \begin{tabular}{  c | c | c | c | c | c | c | c | c }
    \hline
    \text{Slot}             & 1    &  2   &  3    &   4  &   5  &  6    &   7   & 8 \\ \hline \hline
$\underline{E}_1^t$    &2.47    &2.27    &2.18    &1.97    &2.28    &2.66    &3.10    &3.38\\
$\underline{E}_2^t$    &2.57    &1.88    &2.16   &1.56    &1.95    &3.07    &3.44    &3.11\\
$\underline{E}_3^t$    &2.32    &2.43    &1.27    &1.39    &2.14    &1.98    &2.68    &4.04\\
$\underline{E}_4^t$    &2.04    &1.92    &2.33    &2.07    &2.13    &2.36    &3.13   &4.16\\
$\underline{E}_5^t$    &2.11    &1.19    &2.26    &2.19    &1.55    &2.71    &3.37    &2.45\\
$\underline{E}_6^t$    &2.01    &2.29    &2.20    &0.98    &2.43    &3.22    &2.74    &3.93\\
    \hline
     $\alpha^t$              & 0.402   & 0.44    &0.724   &1.32    &1.166  & 0.798   & 0.506   & 0.468\\
    \end{tabular}
\end{table}

\begin{figure}[t]
\centering
\includegraphics[scale=0.5]{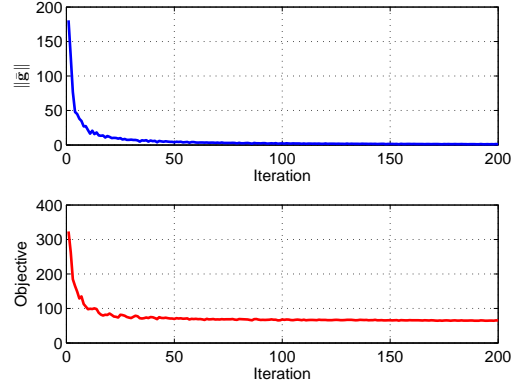}
\caption{Convergence of the objective value and the subgradient norm of Lagrange multipliers
($M=2$, $K=10$, $r=1$).}
\label{fig:S_convg}
\end{figure}

\begin{figure}[t]
\centering
\includegraphics[scale=0.5]{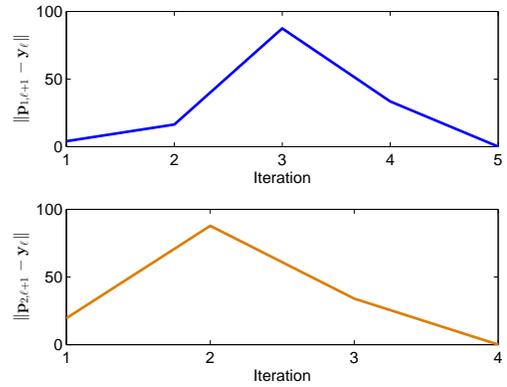}
\caption{Convergence of the bundle method solving the subproblem~\eqref{eq.sub3}
($M=2$, $K=10$, $r=1$, $j=2$).}
\label{fig:S_bundle_convg}
\end{figure}

\begin{figure}[th]
\centering
\includegraphics[scale=0.5]{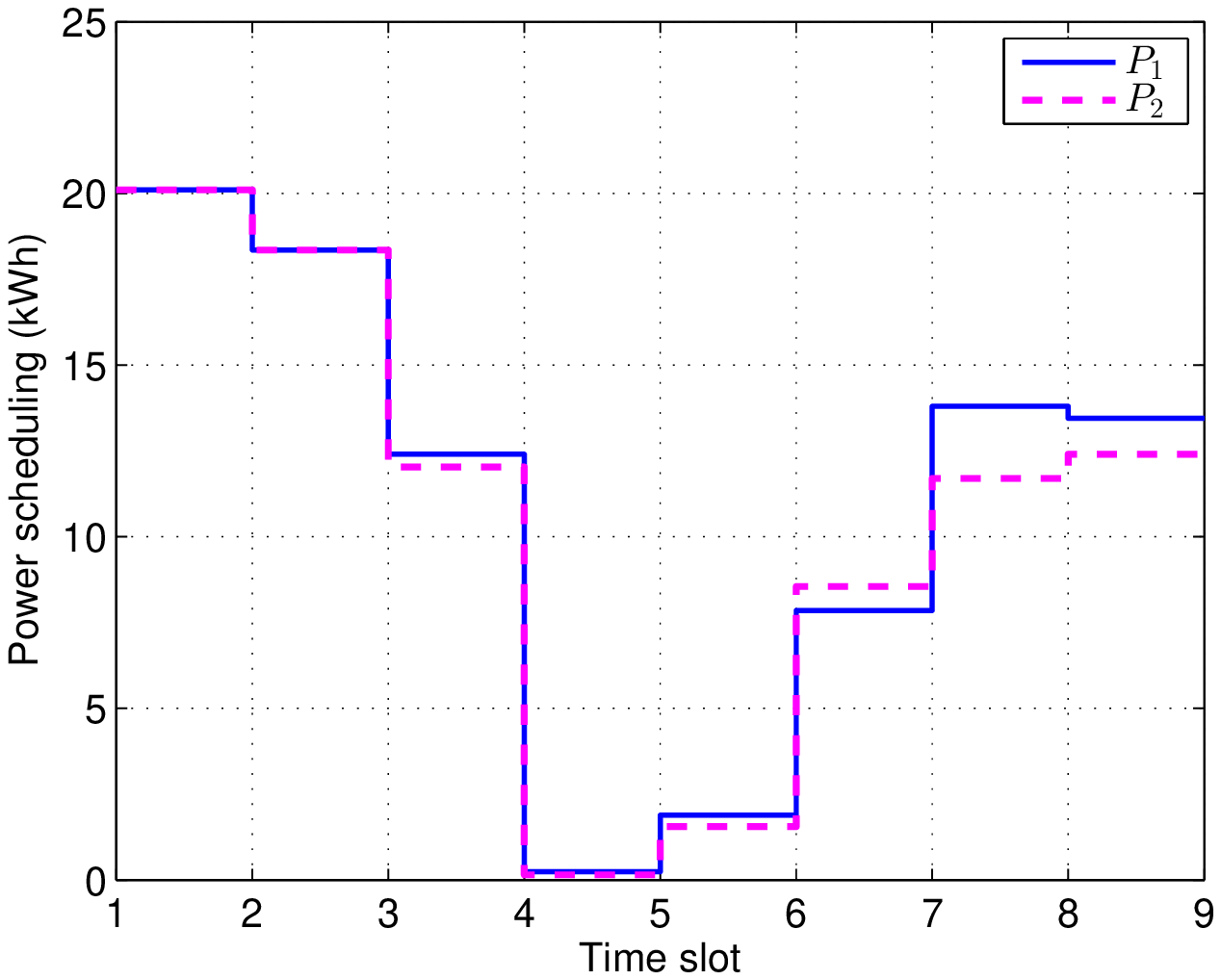}
\caption{Optimal power schedules of $\bar{P}_i^t$ ($I=2$, $M=2$, $K=10$, $r=1$).}
\label{fig:S_pi}
\end{figure}

\begin{figure}[th]
\centering
\includegraphics[scale=0.5]{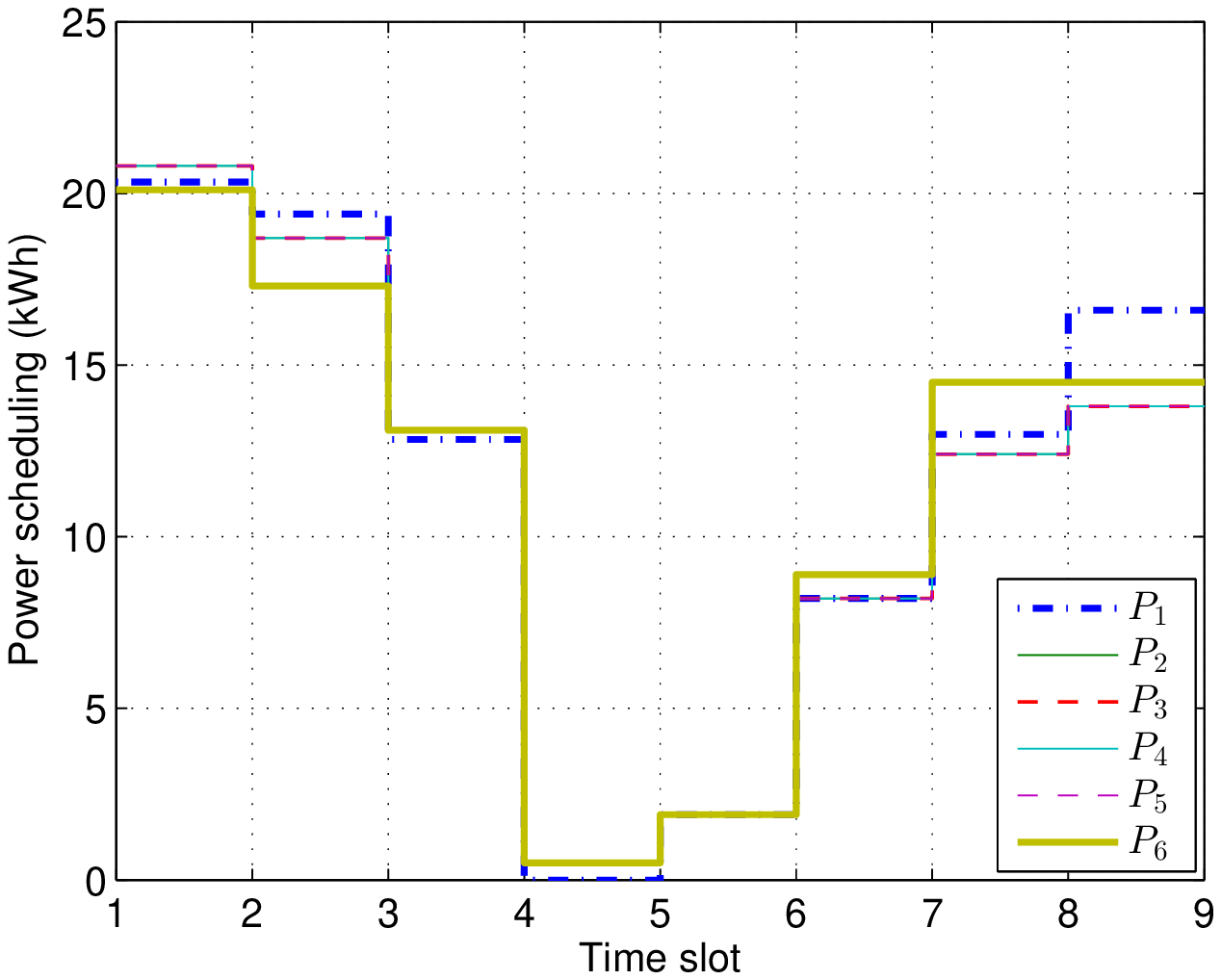}
\caption{Optimal power schedules of $\bar{P}_i^t$ ($I=6$, $M=2$, $K=20$, $r=1$).}
\label{fig:pi6}
\end{figure}

\begin{figure}[th]
\centering
\includegraphics[scale=0.5]{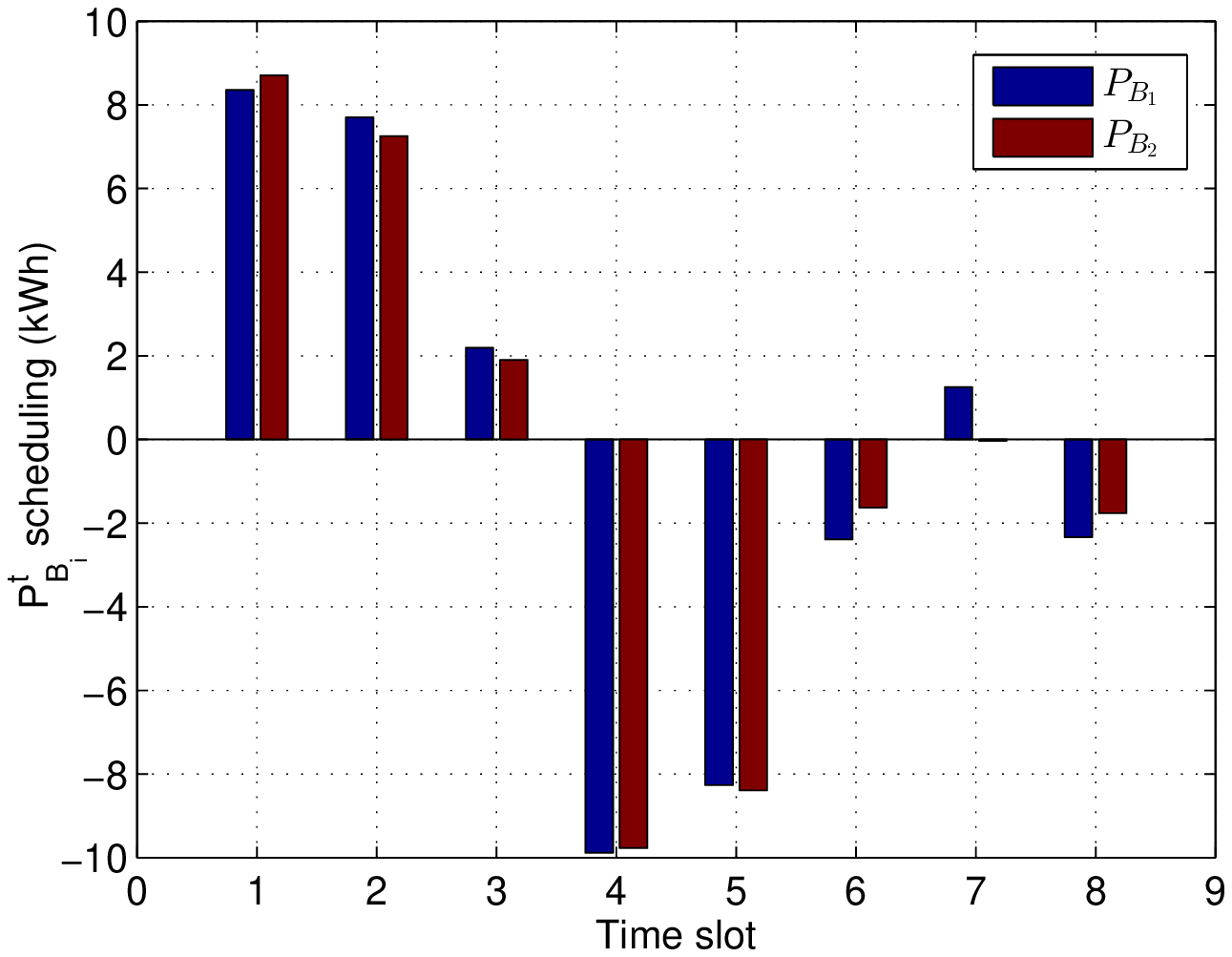}
\caption{Optimal power schedule of $\bar{P}_{B_i}^t$ ($I=2$, $M=2$, $K=10$, $r=1$).}
\label{fig:S_pb}
\end{figure}

\begin{figure}[th]
\centering
\includegraphics[scale=0.5]{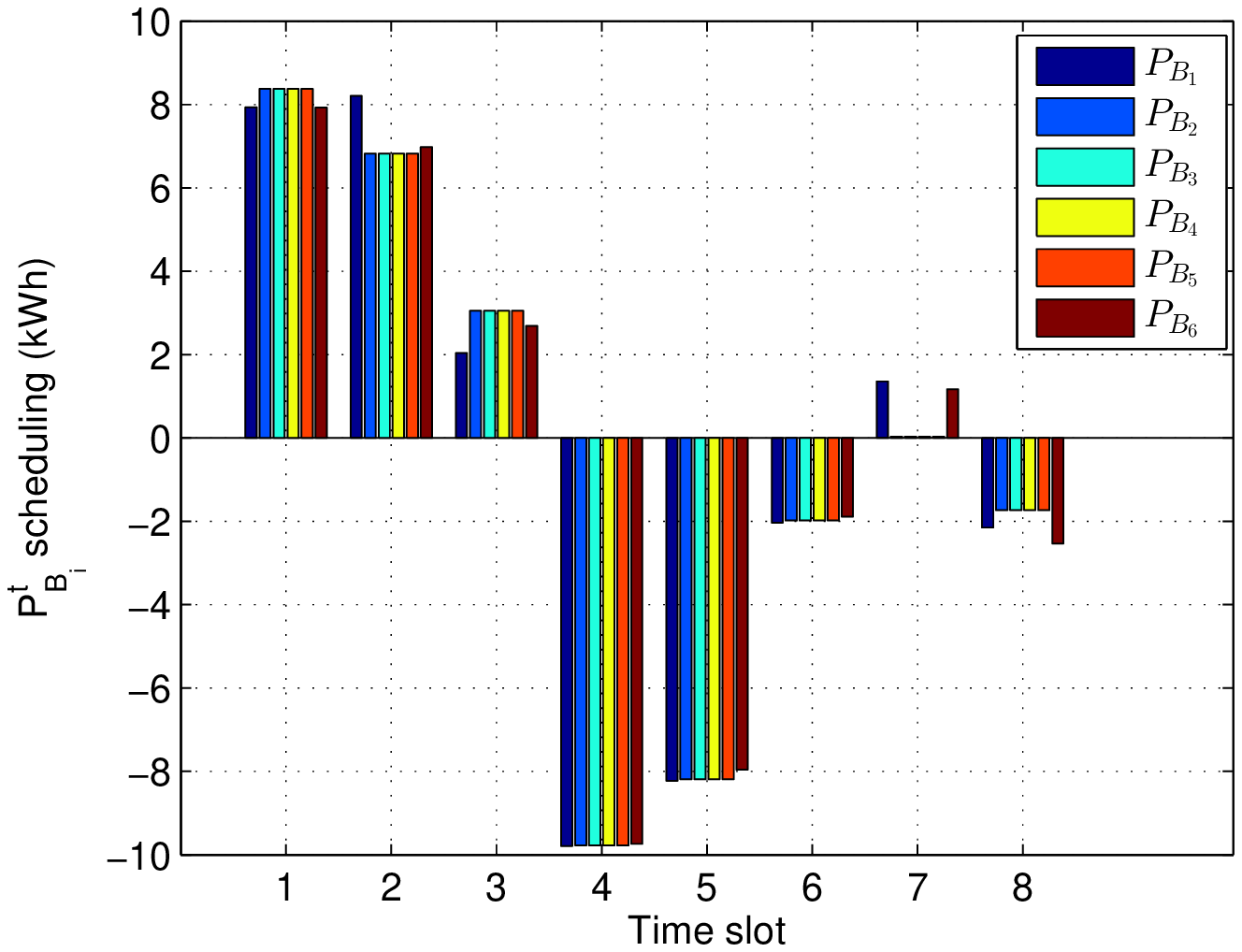}
\caption{Optimal power schedule of $\bar{P}_{B_i}^t$ ($I=6$, $M=2$, $K=20$, $r=1$).}
\label{fig:pb6}
\end{figure}

\begin{figure}[th]
\centering
\includegraphics[scale=0.5]{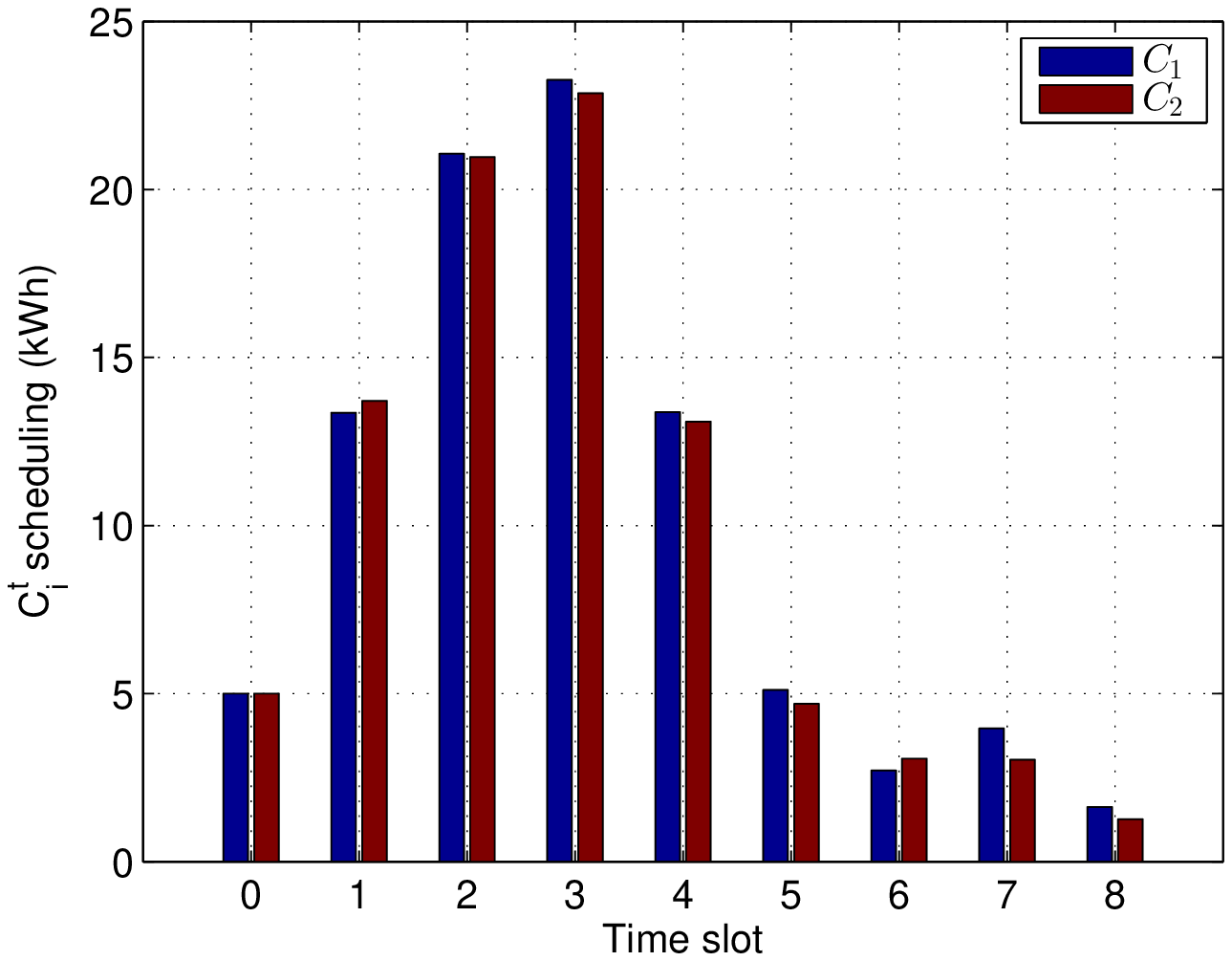}
\caption{Optimal power schedule for $\bar{C}_i^t$ ($I=2$, $M=2$, $K=10$, $r=1$).}
\label{fig:S_soc}
\end{figure}

\begin{figure}[th]
\centering
\includegraphics[scale=0.5]{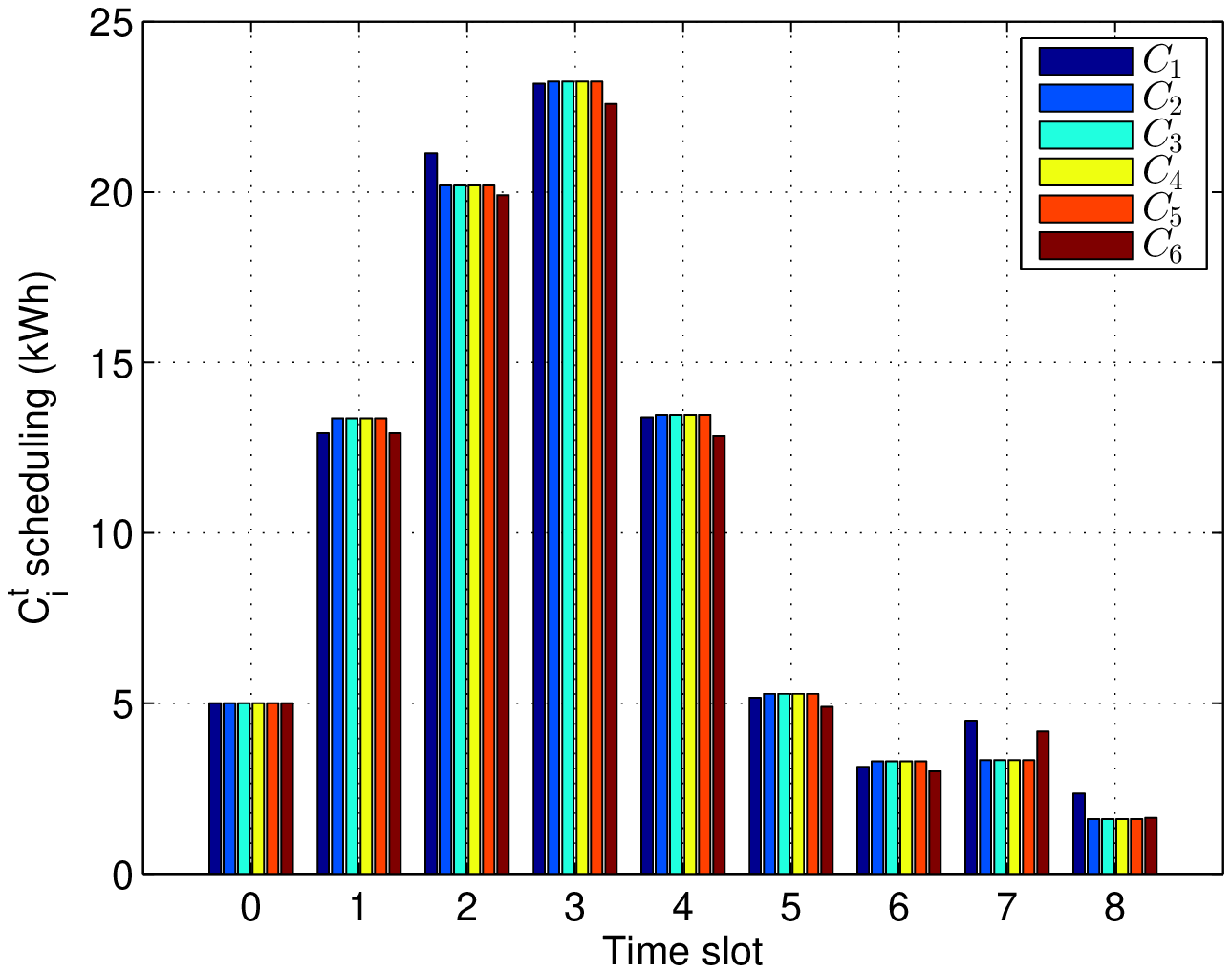}
\caption{Optimal power schedule for $\bar{C}_i^t$ ($I=6$, $M=2$, $K=20$, $r=1$).}
\label{fig:soc6}
\end{figure}

\begin{figure}[th]
\centering
\includegraphics[scale=0.5]{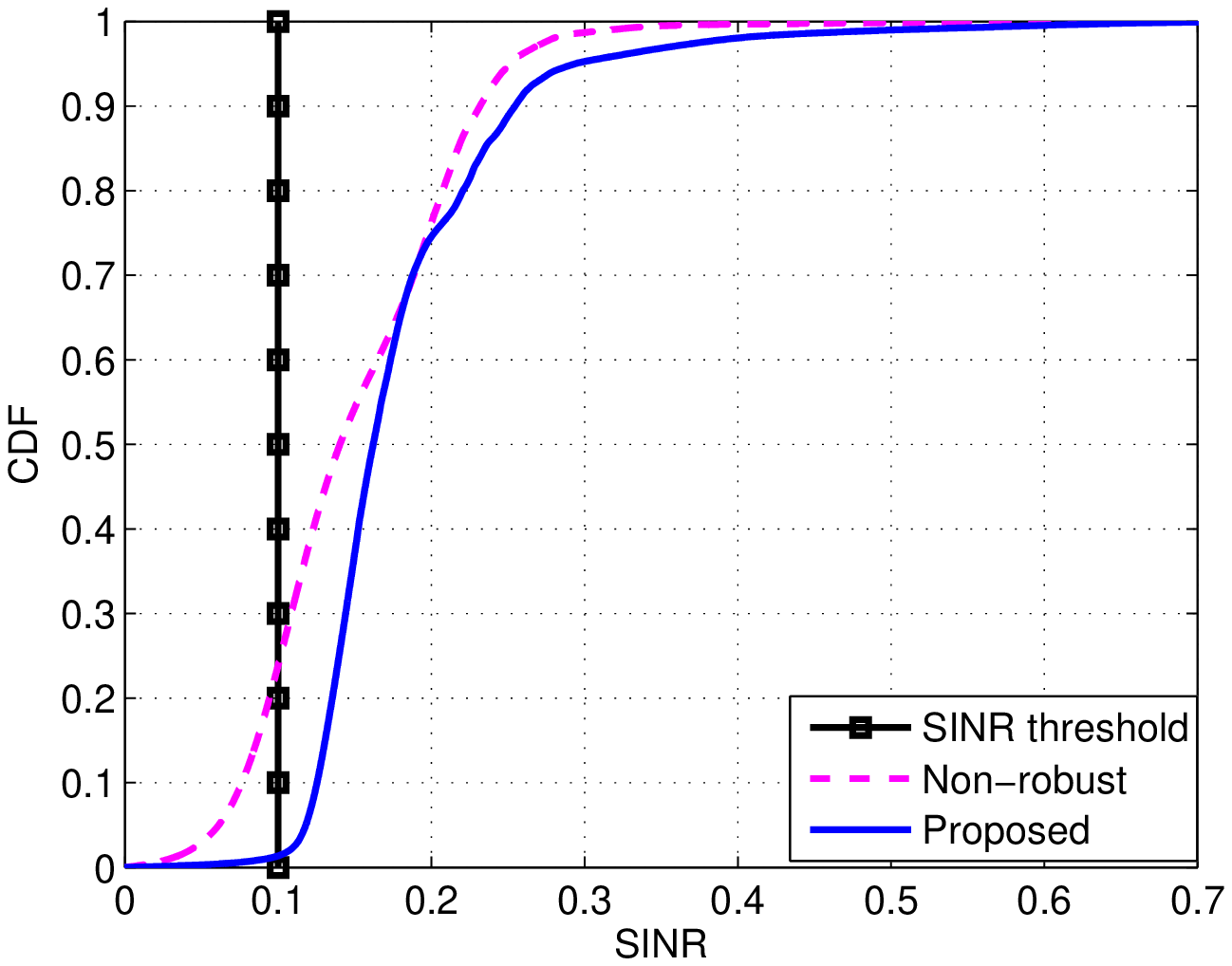}
\caption{SINR cumulative distribution function (CDF) ($I=2$, $M=2$, $K=10$, $r=1$).}
\label{fig:S_sinr_cdf_Xbar}
\end{figure}

%

\begin{figure}[th]
\centering
\includegraphics[scale=0.5]{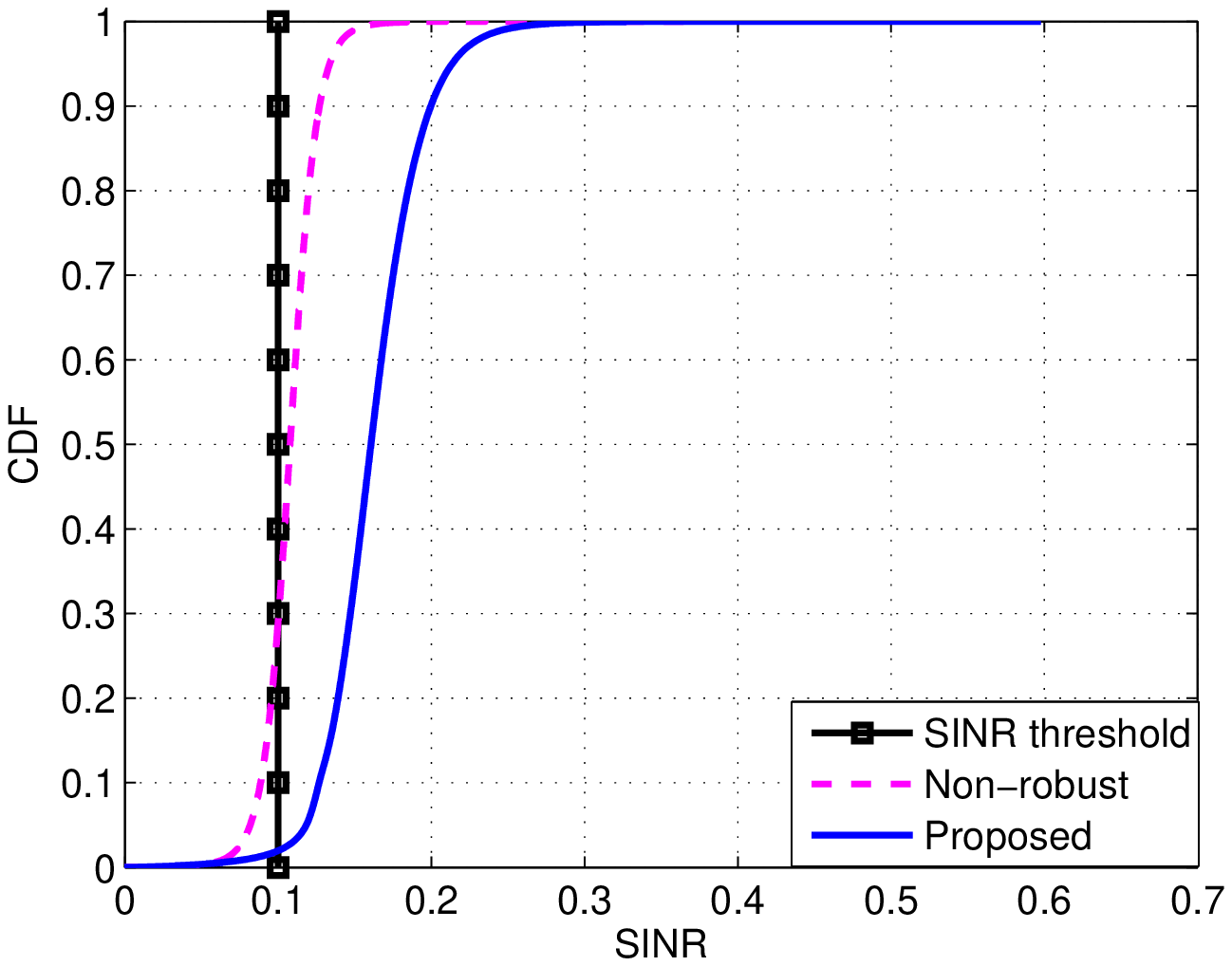}
\caption{SINR cumulative distribution function (CDF) ($I=6$, $M=2$, $K=20$, $r=1$).}
\label{fig:cdf6}
\end{figure}

\begin{figure}[th]
\centering
\includegraphics[scale=0.5]{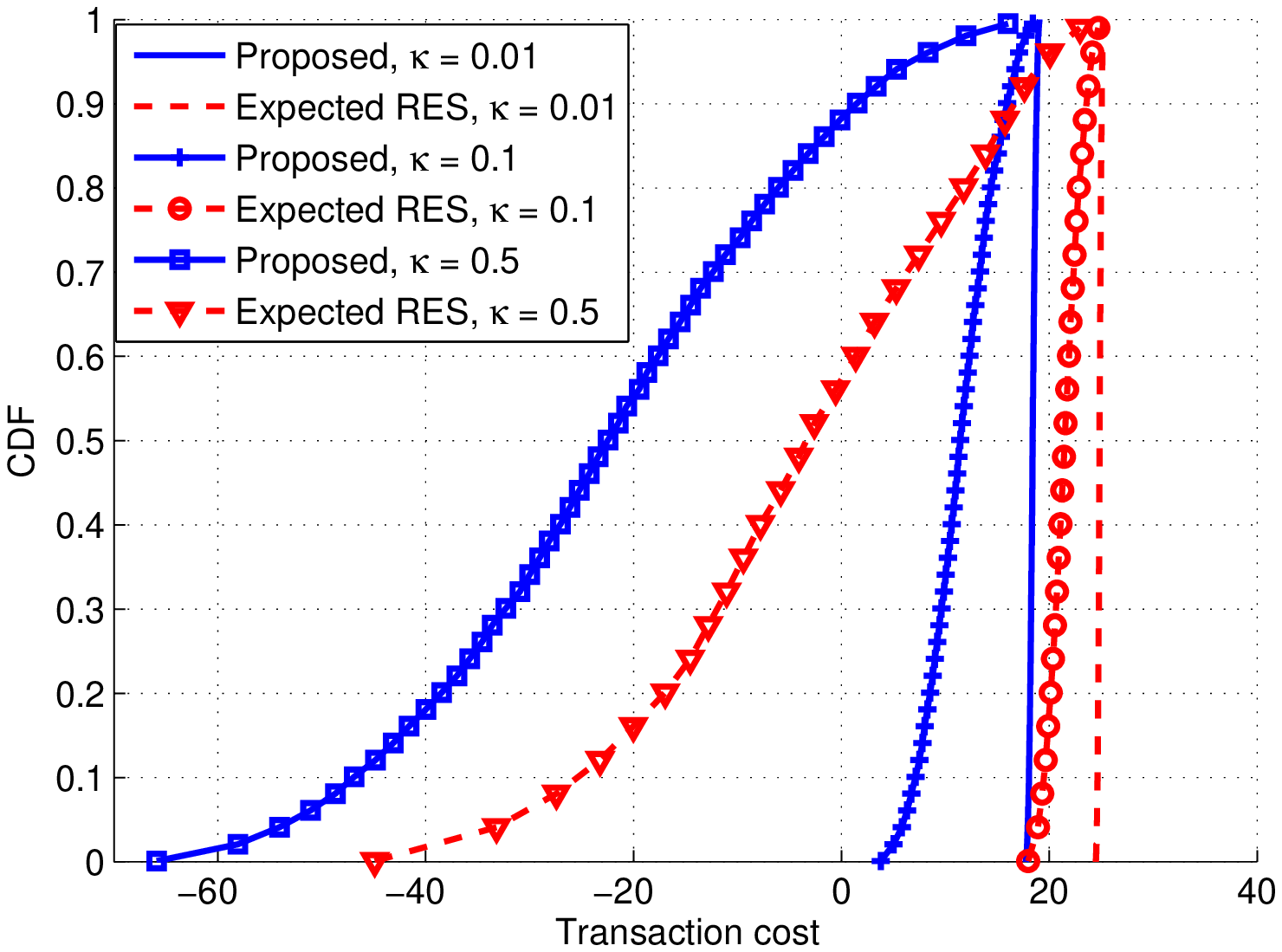}
\caption{Cumulative distribution function (CDF) of the transaction cost ($I=2$, $M=2$, $K=10$, $r=0.3$).}
\label{fig:costComp}
\end{figure}

In this section, simulated tests are presented to verify the performance of the proposed approach.
The Matlab-based modeling package~\texttt{CVX 2.1}~\cite{cvx} along with the solvers~\texttt{MOSEK 7.0}~\cite{mosek}
and~\texttt{Sedumi 1.02}~\cite{SeDuMi} are used to specify and solve the resulting optimization problems.
All numerical tests are implemented on a computer workstation with Intel cores 3.40 GHz and 32 GB RAM.

The scheduling horizons of the considered CoMP network is $T=8$. Two configurations are tested:
(C1) $I=2$ BSs and $K=10$ end users (small size); and (C2) $I=6$ BSs and $K=30$ end users (large size).
All wireless channels are assumed to be flat Rayleigh fading, and normalized to unit power.
The noises are modeled as circularly-symmetric Gaussian random vectors.
Without loss of generality, the effects of path loss, shadowing, and Doppler fading are ignored.
Parameters including the limits of $P_{G_i}$, $C_i$, $P_{b_i}$ and the discharging efficiency
$\varpi_i$ are listed in Table~\ref{tab:param}.
A polyhedral uncertainty set~\eqref{eq.Ei} with a single sub-horizon (no time partition) is considered for the RES.
In Table~\ref{tab:wind}, the energy purchase price $\alpha^t$ is given across the entire time horizon.
The selling price is set as $\beta^t = r\alpha^t$ with $r\in [0,1]$.
In addition, the lower limits $\{\underline{E}_i^t\}_{i\in \mathcal{I}}$ are listed therein,
which were rescaled from real data we obtained from the Midcontinent Independent System Operator (MISO)~\cite{WindData}.
The upper limits were set to $\overline{E}_i^t = 10\underline{E}_i^t$,
while the total horizon bounds are $E_i^{\max}=0.9\sum_{t}\overline{E}_i^t$.

First, convergence of the objective value~\eqref{eq.Reforma}, and the $\ell_2$-norm of the subgradient of
the running-average Lagrange multiplier~\eqref{eq.subgMlpr} is verified in Fig.~\ref{fig:S_convg}.
It can be seen that both metrics converge within a few hundred
iterations, which confirms the validity of the dual decomposition
approach along with the subgradient ascent algorithm. With the Ces\`{a}ro
averages, convergence of the dual and primal variables was
also confirmed, but it is omitted due to limited space.

Figure~\ref{fig:S_bundle_convg} depicts the effectiveness of the proposed
bundle method minimizing the convex nonsmooth objective~\eqref{eq.sub3}.
Clearly, incorporating the scheme of dynamically changing the proximity weight $\rho_{\ell}$,
the bundle algorithm converges very fast; typically, within 10 iterations.

The optimal power schedules of $\bar{P}_i^t$ are depicted in Figs.~\ref{fig:S_pi} and~\ref{fig:pi6}.
For both configurations, the stairstep curves show that the lowest levels of $\bar{P}_i^t$ occur from
slot $4$ to $6$. This is because the energy selling and purchase prices are relatively high
during these horizons [cf.~Sec.~\ref{tab:wind}], which drives the BSs's power consumption low in order to minimize the transaction cost.

Figs.~\ref{fig:S_pb} and~\ref{fig:pb6} show the optimal battery (dis-)charging amount $\bar{P}_{B_i}^t$,
while Figs.~\ref{fig:S_soc} and~\ref{fig:soc6} depict the state of charge $\bar{C}_i^t$.
As a component part of $\bar{P}_i^t$, $\bar{P}_{B_i}^t$ exhibits a similar trend in response
to the price fluctuation; that is, there is a relative large amount of battery discharging ($\bar{P}_{B_i}^t<0$)
when the corresponding price is high.
Note that both $\bar{P}_{B_i}^t$ and $\bar{C}_i^t$ never exceed
their lower and upper limits [cf.~Table \ref{tab:param}].

Robustness of the worst-case design to the uncertainty of channel estimates
[cf.~\eqref{eq.chUncerty}] is confirmed in Figs.~\ref{fig:S_sinr_cdf_Xbar} and~\ref{fig:cdf6}.
The red solid line indicates the SINR threshold $\gamma_k =0.1$ that is set common to all users for simplicity.
The non-robust scheme simply considers the estimated channel $\hat{\mathbf{h}}_k^t$ as the actual one,
and plugs it into the worst-case SINR control design~\eqref{eq.sinr}.
This constraint can be relaxed to a linear matrix inequality:
$\hat{\mathbf{h}}_k^{tH}\mathbf{Y}_k^t\hat{\mathbf{h}}_k^t-\sigma_k^2 \geq 0$,
which is a relaxed version of our proposed counterpart in~\eqref{S-procedure-Cons}.
For both the robust and non-robust approaches, each transmit beamformer $\mathbf{w}_k^t$ is
obtained as the principal eigenvector associated with the largest eigenvalue of the resulting
$\bar{\mathbf{X}}_k^t$.  The CDF of the actual SINR is obtained by evaluating~\eqref{sinr} with $5,000$
independent and identically distributed (i.i.d.) channel realizations.
The channel perturbations $\{\bm{\delta}_k^t\}$ are first generated as complex Gaussian,
and then rescaled to the boundary of the spherical region $\mathcal{H}_k^t$ [cf.~\eqref{eq.chUncerty}].
It can be seen that $20$\% of the realizations of the non-robust scheme does not satisfy the SINR constraint,
while only about $2$\% for the proposed approach. Note that the SDP relaxation is not always exact,
which results in violating the SINR threshold for a few channel realizations.

Finally, CDFs of the transaction cost are depicted in Fig.~\ref{fig:costComp}.
The proposed robust approach is compared with a heuristic scheme that
assumes the expected renewable generation $\hat{E}_i^t = \frac{1}{2}(\underline{E}_i^t+\overline{E}_i^t)$
is available and prefixed for problem~\eqref{eq.sdr}.
The CDF curves were plotted by evaluating the transaction cost~\eqref{eq.Reforma} with $10^5$ realizations
of $\{E_i^t\}_{i,t}$ and the obtained optimal solutions $\{\bar{P}_i^t\}_{i,t}$.
The RES realizations were generated as
$\{\tilde{E}_i^t\}_{i,t} = \underline{E}_i^t + \kappa U(\overline{E}_i^t-\underline{E}_i^t)$,
where $U$ is a uniform random variable on $[0, 1]$.
Three cases with $\kappa = 0.01, 0.1, 0.5$ were tested.
Clearly, transaction costs of both the proposed and the heuristic methods decrease with the increase of $\kappa$.
Since a larger value of $\kappa$ implies more harvested renewables energy yielding a reduced transaction cost.
Note that negative transaction costs means net profits are obtained by selling surplus renewables back to the smart grid.
Interestingly, for all cases, the proposed approach always outperforms the heuristic scheme with less transaction costs.
This is because the proposed schedules of $\{P_{g,i}^t,P_{b,i}^t\}_{i,t}$ are robust to the worst-case renewable generation
$\{E_i^t\}_{i,t}$. However, in practice more RES is often available than the worst case.
Hence, the proposed method has a larger profit-making capability than the heuristic alternative.

\section{Conclusions}

Robust ahead-of-time energy management and transmit-beamforming designs were developed to minimize the worst-case energy costs subject to the worst-case user QoS guarantees for the CoMP downlink, which is powered by a grid with smart-meter based dynamic pricing and RES available at the BSs. Building on innovative models, the task was formulated as a convex program. Relying on a Lagrange-dual approach, efficient algorithms were introduced to obtain optimal solutions. The proposed scheme provides the offline ahead-of-time beamformer and energy schedule over a finite time horizon. Development of online management schemes will be an interesting direction to pursue in our future work.

Supported by major programs such as the US Energy Independence and Security Act and European SmartGrids Technology Platform, the smart grid industry has seen fast growth. It is expected that all future communication systems will be powered by a smart grid. As a result, integration of smart-grid technologies into system designs will hold the key to fully leverage energy-efficient communications in next-generation HetNets. The proposed models and approaches can pave a way to further advancing fundamental research on smart-grid powered HetNet transmissions, which will be pursued in future works.

\bibliographystyle{IEEEtran}
\bibliography{gm15_bib}

\end{document}